\documentclass{amsart}
\usepackage[margin=2cm,letterpaper]{geometry}
\usepackage{amsaddr,amsthm,amssymb,enumerate,stackrel}
\usepackage{xr-hyper}
\newcommand{\csp}[1]{\begin{NoHyper}\cite[#1]{PavlovScholbach:Spectra}\end{NoHyper}}
\newcommand{\cop}[1]{\begin{NoHyper}\cite[#1]{PavlovScholbach:Operads}\end{NoHyper}}
\usepackage[hypertex,backref=section]{hyperref}
\usepackage[all,cmtip]{xy}
\SilentMatrices 
\usepackage[utf8]{inputenc} 

\def\typeout#1{} 

\makeatletter
\@addtoreset{equation}{section}

\makeatother

\theoremstyle{definition}
\newtheorem{Defi}[equation]{Definition} \newcommand{\defi}{\begin{Defi}} \newcommand{\xdefi}{\end{Defi}} \newcommand{\refde}[1]{Definition~\ref{defi--#1}}
\newtheorem{Bsp}[equation]{Example} \newcommand{\exam}{\begin{Bsp}} \newcommand{\xexam}{\end{Bsp}} 

\theoremstyle{remark}
\newtheorem{Bem}[equation]{Remark} \newcommand{\rema}{\begin{Bem}} \newcommand{\xrema}{\end{Bem}} \newcommand{\refre}[1]{Remark~\ref{rema--#1}}
\newtheorem{Nota}[equation]{Notation} \newcommand{\nota}{\begin{Nota}} \newcommand{\xnota}{\end{Nota}} 

\theoremstyle{plain}
\newtheorem{Theo}[equation]{Theorem} \newcommand{\theo}{\begin{Theo}} \newcommand{\xtheo}{\end{Theo}} \newcommand{\refth}[1]{Theorem~\ref{theo--#1}}
\newtheorem{Satz}[equation]{Proposition} \newcommand{\prop}{\begin{Satz}} \newcommand{\xprop}{\end{Satz}} \newcommand{\refpr}[1]{Proposition~\ref{prop--#1}}
\newtheorem{Lemm}[equation]{Lemma} \newcommand{\lemm}{\begin{Lemm}} \newcommand{\xlemm}{\end{Lemm}} \newcommand{\refle}[1]{Lemma~\ref{lemm--#1}}
\newtheorem{Coro}[equation]{Corollary} \newcommand{\coro}{\begin{Coro}} \newcommand{\xcoro}{\end{Coro}} \newcommand{\refco}[1]{Corollary~\ref{coro--#1}}

\newcommand{\refsect}[1]{\S\ref{sect--#1}}
\newcommand{\refchap}[1]{Section~\ref{chap--#1}}
\newcommand{\refit}[1]{(\ref{item--#1})}
\newcommand{\refeq}[1]{(\ref{eqn--#1})}

\newcommand{\pf}{\begin{proof}} \newcommand{\xpf}{\end{proof}}

\let\ppar\endgraf 

\newcommand{\category}[1]{\mathbf{#1}}
\newcommand{\Set}{\category{Set}} 
\newcommand{\simpl}{\category s} 
\newcommand{\sSet}{\simpl\category{Set}} 
\newcommand{\PSh}{\category{PSh}} 
\newcommand{\Top}{\category{Top}} 
\newcommand{\sPSh}{\simpl\category{PSh}} 
\newcommand{\sAb}{\simpl\category{Ab}} 
\newcommand{\sMod}{\simpl\category{Mod}} 
\newcommand{\Mod}{\category{Mod}} 
\newcommand{\Alg}{\mathbf{Alg}} 
\newcommand{\Ch}{\category{Ch}} 
\newcommand{\Ar}{\category{Ar}} 
\newcommand{\Fun}{\category{Fun}} 
\newcommand{\Orb}{\category{Orb}} 

\newcommand{\proj}{\mathrm{pro}} 
\newcommand{\loc}{\mathrm{loc}} 
\newcommand{\inje}{\mathrm{in}} 

\newcommand{\cof}{\mathrm{C}} 
\newcommand{\cofh}{\cof_{\mathrm h}}
\newcommand{\cofi}{\cof_{\mathrm i}}
\newcommand{\dom}{\mathrm{dom}} 
\newcommand{\fib}{\mathrm{F}} 
\newcommand{\we}{\mathrm{W}} 
\newcommand{\AF}{\mathrm{AF}} 
\newcommand{\AC}{\mathrm{AC}} 
\newcommand{\Acofh}{\AC_{\mathrm h}}
\newcommand{\Acofi}{\AC_{\mathrm i}}
\newcommand{\cell}{\mathrm{cell}} 
\def\mcr{{\rm Q}} 
\def\rdf{{\underline{\rm R}}} 

\newcommand{\ws}[1]{\mathop{\rm cof}({#1})} 
\newcommand{\inj}[1]{\mathop{\rm inj}({#1})} 
\newcommand{\colim}{\operatornamewithlimits{colim}} 
\newcommand{\hocolim}{\operatornamewithlimits{hocolim}} 
\newcommand{\coeq}{\operatorname{coeq}} 
\def\id{{\rm id}} 
\def\ev{{\rm ev}} 
\def\op{{\rm op}} 
\def\Mor{\mathop{\rm Mor}\nolimits} 
\def\Hom{\mathop{\rm Hom}\nolimits} 
\def\hHom{\mathop{\rdf \rm Hom}\nolimits} 
\def\hMap{\mathop{\rdf \rm Map}\nolimits} 
\def\Ax{\Sigma} 

\mathchardef\pp"2403 
\mathchardef\ppdom"2400 

\def\bigpp{\mathop{\mathchoice{\dobigpp\Huge}{\dobigpp\Large}{\dobigpp\normalsize}{\dobigpp\small}}}
\mathchardef\bigppchar"1403
\def\dobigpp#1{\vcenter{#1\kern.2ex\hbox{$\bigppchar$}\kern.2ex}}
\mathchardef\bigppdomchar"1400
\def\bigppdom{\bigppdomchar\nolimits} 
\def\dobigppdom#1{\vcenter{#1\kern.2ex\hbox{$\bigppdomchar$}\kern.2ex}}
\def\bigppdom{\mathop{\mathchoice{\dobigppdom\Huge}{\dobigppdom\Large}{\dobigppdom\normalsize}{\dobigppdom\small}}\nolimits}

\def\Z{{\bf Z}} 
\def\DD{{\mathrm D}} 
\def\Q{{\bf Q}} 
\def\RR{{\bf R}} 
\def\P{{\bf P}} 

\def\H{{\rm H}} 
\def\BL{{\rm L}} 
\def\LL{{\rm L}} 

\newcommand{\Comm}{\mathrm{Comm}} 
\def\Ei{{\rm E_\infty}} 
\def\E{{\rm E}} 
\def\DC{{\rm DC}} 
\newcommand{\Y}{\mathcal{Y}} 
\newcommand{\C}{\mathcal{C}} 
\newcommand{\D}{\mathcal{D}} 
\def\V{{\cal V}} 
\newcommand{\cA}{\mathcal{A}} 
\newcommand{\cO}{\mathcal{O}}

\def\To#1#2{\mathop{\count0=#1 \loop\ifnum\count0>0 \smash-\mkern-7mu \advance\count0 -1 \repeat \mathord\rightarrow}\limits^{#2}} 
\let\x\times

\let\t\otimes
\let\r\rightarrow

\def\matrix#1{\null\,\vcenter{\normalbaselines
    \ialign{\hfil$##$\hfil&&\quad\hfil$##$\hfil\crcr
      \mathstrut\crcr\noalign{\kern-\baselineskip}
      #1\crcr\mathstrut\crcr\noalign{\kern-\baselineskip}}}\,}
\def\vcd#1{\def\normalbaselines{\baselineskip20pt\lineskip1pt\lineskiplimit0pt } \harrowsize#1 \matrix}
\def\cd{\vcd{2em}}
\newdimen\harrowsize
\def\mapright#1{\smash{\mathop{\hbox to\harrowsize{\rightarrowfill}}\limits^{#1}}}

\def\mapdown#1{\Big\downarrow\rlap{$\vcenter{\hbox{$\scriptstyle#1$}}$}}

\catcode`@=11
\let\over\@@over
\let\atop\@@atop
\let\above\@@above
\let\overwithdelims\@@overwithdelims
\let\atopwithdelims\@@atopwithdelims
\let\abovewithdelims\@@abovewithdelims
\def\eqalign#1{\null\,\vcenter{\openup\jot\m@th
  \ialign{\strut\hfil$\displaystyle{##}$&$\displaystyle{{}##}$\hfil
      \crcr#1\crcr}}\,}
\newskip\xcentering
\def\eqalignno#1{\displ@y \tabskip\xcentering
  \halign to\displaywidth{\hfil$\@lign\displaystyle{##}$\tabskip\z@skip
    &$\@lign\displaystyle{{}##}$\hfil\tabskip\xcentering
    &\llap{$\@lign##$}\tabskip\z@skip\crcr
    #1\crcr}}
\def\cases#1{\left\{\,\vcenter{\normalbaselines\m@th
    \ialign{$##\hfil$&\quad##\hfil\crcr#1\crcr}}\right.}
\def\@writetocindents{}
\def\eqlabel#1{\refstepcounter{equation}\label{eqn--#1}\ifmmode\ifinner\else\eqno\fi\fi\hbox{\@eqnnum}} 
\def\cal{\fam\tw@}

\def\uppercasenonmath#1{}

\catcode`@=12

\begin{document}

\title{Homotopy theory of symmetric powers}

\author{Dmitri Pavlov}
\address{Faculty of Mathematics, University of Regensburg\\
Department of Mathematics and Statistics, Texas Tech University}

\author{Jakob Scholbach}
\address{Mathematical Institute, University of M\"unster}


\centerline{\font\tfont=cmss17 \tfont\shorttitle}

\bigskip
{\tabskip0pt plus 1fil
\halign to\hsize{&\hfil#\hfil\cr
{\bf Dmitri Pavlov\/}\cr
Faculty of Mathematics, University of Regensburg\cr
Department of Mathematics and Statistics, Texas Tech University\cr
\href{https://dmitripavlov.org/}{https:/\negthinspace/dmitripavlov.org/}\vadjust{\medbreak}\cr
{\bf Jakob Scholbach\/}\cr
Mathematical Institute, University of M\"unster\cr
\href{https://wwwmath.uni-muenster.de/u/jakob.scholbach/}{https:/\negthinspace/wwwmath.uni-muenster.de/u/jakob.scholbach/}\cr
}}

\begin{abstract}
We introduce the symmetricity notions of symmetric h-monoidality, symmetroidality, and symmetric flatness.
As shown in our paper \cite{PavlovScholbach:Operads}, these properties lie at the heart of the homotopy theory of colored symmetric operads and their algebras.
In particular, the former property can be seen as the analog of Schwede and Shipley's monoid axiom for algebras over symmetric operads
and allows one to equip categories of such algebras with model structures,
whereas the latter ensures that weak equivalences of operads induce Quillen equivalences of categories of algebras.
We discuss these properties for elementary model categories such as simplicial sets, simplicial presheaves, and chain complexes.
Moreover, we provide powerful tools to promote these properties from such basic model categories to more involved ones, such as the stable model structure on symmetric spectra.
\end{abstract}


\makeatletter\@setabstract\makeatother

\setcounter{tocdepth}{2}
\tableofcontents

\numberwithin{equation}{section}

\section{Introduction}

Model categories provide an important framework for homotopy-theoretic computations.
Algebraic structures such as monoids, their modules, and more generally operads and their algebras provide means
to concisely encode multiplication maps and their properties such as unitality, associativity, and commutativity.
Homotopy coherent versions of such algebraic structures form the foundation of a variety of mathematical areas,
such as algebraic topology, homological algebra, derived algebraic geometry, higher category theory, and derived differential geometry.
This motivates the following question:
what conditions on a monoidal model category~$(\C,\otimes)$ are needed for a meaningful homotopy theory of monoids, modules, etc.?
The first answer to this type of question was given by Schwede and Shipley's \emph{monoid axiom},
which guarantees that for a monoid~$R$ in~$\C$, the category $\Mod_R(\C)$ of $R$-modules carries a model structure transferred from~$\C$, see~\cite{SchwedeShipley:Algebras}.
The monoid axiom asks that transfinite compositions of pushouts of maps of the form $Y\t s$,
where $s$~is an acyclic cofibration and $Y$~is any object are again weak equivalences.
Moreover, given two weakly equivalent monoids $R\buildrel\sim\over\to S$, the categories $\Mod_R$ and $\Mod_S$ are Quillen equivalent
if $Y \t X \r Y' \t X$ is a weak equivalence for any weak equivalence $Y \r Y'$ and any cofibrant object~$X$.

This paper is devoted to a thorough study of the homotopy-theoretic behavior of more general algebraic expressions in a model category, such as
$$X^{\t n}_{\Sigma_n},\qquad Y\otimes_{\Sigma_n}X^{\otimes n},\qquad Z\otimes_{\Sigma_{n_1} \x \cdots \x \Sigma_{n_e}} (X_1^{\t n_1} \t \cdots \t X_e^{\t n_e}),\eqlabel{XYZ}$$
where $X, Y, Z \in \C$, $Y$~has an action of~$\Sigma_n$, $Z$~has an action of~$\prod \Sigma_{n_i}$, and the subscripts denote coinvariants by the corresponding group actions.
More specifically, we introduce \emph{symmetricity properties} for a symmetric monoidal model category $\C$:
\emph{symmetric h-monoidality}, \emph{symmetroidality}, and \emph{symmetric flatness}.

Symmetric h-monoidality requires, in particular, that for any object~$Y$ as above and any acyclic cofibration~$s$ in~$\C$, the map
$$Y\otimes_{\Sigma_n}s^{\pp n}\eqlabel{symmetric.h.monoidal.intro}$$
is a couniversal weak equivalence, i.e., a map whose cobase changes are weak equivalences.
Here $s^{\pp n}$ is the $n$-fold pushout product of~$s$, which is a monoidal product on morphisms.
Symmetric h-monoidality is a natural enhancement of h-monoidality introduced by Batanin and Berger in~\cite{BataninBerger:Homotopy}.

Symmetric flatness requires that for any $\Sigma_n$-equivariant map~$y$ whose underlying map in~$\C$ is a weak equivalence and any cofibration $s \in \C$, the map
$$y \pp_{\Sigma_n}s^{\pp n}\eqlabel{symmetric.flat.intro}$$
is a weak equivalence.
This implies that $y \t_{\Sigma_n} X^{\t n}$ is a weak equivalence for any cofibrant object~$X$.
Among other things this means that the $\Sigma_n$-quotients in~\refeq{XYZ} are also homotopy quotients.
See \ref{defi--symmetric.i.monoidal},~\ref{defi--symmetric.flat} for the precise definitions.

Expressions as in~\refeq{XYZ} are of paramount importance for handling monoids and, more generally, algebras over colored symmetric operads.
Indeed, a free commutative monoid, more generally, a free algebra over a (colored) symmetric operad, involves such terms.
In \cite{PavlovScholbach:Operads}, we show that symmetric h-monoidality ensures the existence of a transferred model structure on algebras over
\emph{any} symmetric colored operad, while symmetric flatness yields a Quillen equivalence of algebras over weakly equivalent operads.
We also introduce \emph{symmetroidality} in this paper, which can be used to govern the behavior of cofibrant algebras over operads.

Up to transfinite compositions present in the monoid axiom, which we treat separately, symmetric h-mon\-oidality and symmetric flatness
can be regarded as natural enhancements of the above conditions of Schwede and Shipley.
However, it turns out to be hard to establish the symmetric h-monoidality, symmetroidality, and symmetric flatness for a given model category $\C$ directly.
Therefore, the main goal of this paper is to provide a powerful and convenient set of tools that enable us to quickly promote these properties through various constructions on model categories.

\theo (See \refth{symmetric.weakly.saturated} for the precise statement.)
To check that $\C$ is symmetric h-monoidal or symmetric flat it is enough to consider~\refeq{symmetric.h.monoidal.intro} and~\refeq{symmetric.flat.intro}
for \emph{generating} cofibrations~$s$.
\xtheo

\theo (See \refth{transfer.symmetric.monoidal} for the precise statement.)
Suppose $F : \C \rightleftarrows \D : G$ is an adjunction of symmetric monoidal model categories
that is sufficiently compatible with the monoidal products, such as $\D = \Mod_R (\C)$, where $R$~is a commutative monoid in~$\C$.
Then the symmetric h-monoidality and symmetric flatness of~$\C$ imply the one of~$\D$.
\xtheo

\theo (See \refth{symmetric.monoidal.localization} for the precise statement.)
Given a monoidal left Bousfield localization
$$\C \rightleftarrows \D = \LL^\t_S (\C),$$
the symmetric h-monoidality and symmetric flatness of $\C$ imply the one of $\D$.
\xtheo

As an illustration of these principles, consider the problem of establishing the symmetric h-monoidality, symmetroidality,
and symmetric flatness for the monoidal model category of simplicial symmetric spectra.
This allows one to establish the homotopy theory of operads and their algebras in spectra, such as commutative ring spectra or $\Ei$-ring spectra.
First, by direct inspection~(\refsect{simplicial.sets}) one establishes these properties for the generating (acyclic) cofibrations of simplicial sets, i.e., $\partial\Delta^n\to\Delta^n$ and $\Lambda^n_k\to\Delta^n$.
By \refth{symmetric.weakly.saturated}, this shows that $\sSet$ is symmetric h-monoidal, symmetroidal, and flat.
Next, again by direct inspection, one can show that positive cofibrations of symmetric sequences (i.e., cofibrations that are isomorphisms in degree~0)
form a symmetric h-monoidal, symmetric flat class.
Via \refth{transfer.symmetric.monoidal} these properties can be transferred to modules over a (fixed) commutative monoid in symmetric sequences
(specifically, the sphere spectrum),
equipped with the positive unstable (i.e., transferred) model structure.
Finally, by applying \refth{symmetric.monoidal.localization}, one establishes them for the left Bousfield localization of the positive unstable model structure with respect to the stabilizing maps,
which gives the positive stable model structure on simplicial symmetric spectra.
These steps are carried out in detail for spectra in an abstract model category in \cite{PavlovScholbach:Spectra}.

We begin by recalling some basic notions pertaining to model categories in \refchap{model} and the monoidal structure on the arrow category $\Ar(\C)$~(\refsect{arrows}).
We then recall the notion of h-monoidality due to Batanin and Berger~\cite{BataninBerger:Homotopy},
and the concept of flatness, which is well-known and has been independently studied by Hovey, for example, see~\cite{Hovey:Smith}.
In \refchap{symmetric}, we define the above-mentioned symmetricity concepts.
This extends the work of Lurie \cite{Lurie:HA} and Gorchinskiy and Guletski\u\i~\cite{GorchinskiyGuletskii:Symmetric}.
An important technical key is \refth{symmetric.weakly.saturated}, which shows the stability of these properties under weak saturation.
This extends a similar statement of Gorchinskiy and Guletski\u\i~\cite[Theorem~5]{GorchinskiyGuletskii:Symmetric}
about stability under weak saturation of a special case of symmetroidality (which we also prove in \refth{symmetric.weakly.saturated}).
More recent accounts of this result include White~\cite[Appendix~A]{White:Model} and Pereira~\cite[\S4.2]{Pereira:Cofibrancy}.
Our proof uses similar ideas, but is shorter.
The stability of the symmetricity and various other model-theoretic properties under transfers and left Bousfield localizations is shown in \S\ref{chap--transfer} and~\S\ref{chap--localization}.
Given that these two methods are the most commonly used tools to construct model structures, the main results of these sections (\ref{prop--transfer.monoidal},
\ref{theo--transfer.symmetric.monoidal}, \ref{prop--monoidal.localization}, \ref{theo--symmetric.monoidal.localization})
should be useful to establish the symmetricity for many other model categories not considered in this paper.
For example, the combination of h-monoidality and flatness allows to carry through the monoid axiom to a left Bousfield localization.
This is illustrated in \refchap{examples}, where we discuss the symmetricity properties of model categories such as simplicial sets, simplicial modules, and simplicial (pre)sheaves, as well as topological spaces and chain complexes.

We thank John Harper, Jacob Lurie, Birgit Richter, Brooke Shipley, and David White for helpful conversations.
This work was partially supported by the SFB 878 grant.

\section{Model categories}
\label{chap--model}

In this section we recall parts of the language of model categories that is used throughout this paper.
As for the definition of a model category we follow \cite{Hovey:Model} or \cite{Hirschhorn:Model}.
The classes of acyclic (co)fibrations are denoted $\AF$ ($\AC$, respectively).
The weak saturation (closure under pushouts, transfinite compositions, and retracts \cite[Definition A.1.2.2]{Lurie:HTT}) of some class~$M$ of morphisms is denoted~$\ws M$.
The class of maps having the right lifting property with respect to all maps in~$M$ is denoted $\inj M$.
Different model structures on the same category are distinguished using superscripts.

\defi \label{defi--model.category.basic}
A model category is
{\it quasi-tractable\/} if its (acyclic) cofibrations are contained in the weak saturation of (acyclic) cofibrations with cofibrant source (and target).

\label{defi--pretty.small}
A model category~$\C$ is \emph{pretty small} if there is a cofibrantly generated model category structure~$\C'$ on the same category as~$\C$ such that $\we_\C=\we_{\C'}$,
$\cof_\C\supset\cof_{\C'}$ and the domains and codomains~$X$ of some set of generating cofibrations of~$\C'$ are compact (also known as finitely presentable), i.e., $\Mor(X,-)$ preserves arbitrary filtered colimits.
\xdefi

Pretty smallness differs from the notion of {\it compact generatedness\/} \cite[Definition~15.2.1]{MayPonto:More} only at the level of technical detail.
It is a convenient assumption ensuring that weak equivalences are stable under colimits of chains (\refle{sequential}).
Moreover, all the basic model categories in \refchap{examples} (except for topological spaces,
which can be treated by a more narrowly tailored compactness condition).
Moreover, pretty smallness is stable under transfer and localization.
Any other compactness condition on a model category satisfying these properties can be used instead of pretty smallness.

\lemm \label{lemm--sequential}
Let $\lambda$ be an ordinal and $f\colon\lambda\to\Ar(\C)$ a cocontinuous chain of morphisms in a model category,
i.e., a sequence of commutative squares
$$\xymatrix @-1pc{X_i \ar[d]_{f_i} \ar[rr]^{x_i} & & X_{i+1} \ar[d]^{f_{i+1}}\cr
Y_i \ar[rr]_{y_i} & & Y_{i+1}\cr}$$ indexed by $i\in\lambda$ such that $f_i=\colim_{j<i}f_j$ for all limit ordinals $i\in\lambda$.
Set $f_\infty=\colim f_i$.
\begin{enumerate}[(i)]
\item \label{item--cof.sequential} \cite[Proposition~I.2.6.3]{ChacholskiScherer:Homotopy}
If every $f_i$ (equivalently, only $f_0$) and every map $X_{i+1}\sqcup_{X_i}Y_i\to Y_{i+1}$ is an (acyclic) cofibration, then so is $f_\infty$.
\item
\label{item--af.sequential}
If cofibrations in~$\C$ are generated by cofibrations with compact domain and codomain
and every~$f_i$ is an acyclic fibration, then so is $f_\infty$.
\item
\label{item--we.chain.colimits}
If $\C$ is pretty small and every~$f_i$ is a weak equivalence, then so is $f_\infty$.
In particular, colimits of chains are homotopy colimits.
The same is true for arbitrary filtered colimits.
\item \label{item--we.transfinite.composition} If $\C$ is pretty small then weak equivalences are stable under transfinite compositions,
i.e., for any cocontinuous chain $X\colon\lambda\to\C$ of weak equivalences the map $X_0\to\colim X$ is also a weak equivalence.
\end{enumerate}
\xlemm
\pf

\refit{af.sequential}: Following the proof of \cite[Corollary~7.4.2]{Hovey:Model}, consider 
$$\xymatrix @-1pc{A \ar[r] \ar[d] & X_s \ar[d] \\
B \ar[r] & Y_s,}$$
where $A\to B$ is a generating cofibration and $s=\infty$.
Since $A$ and $B$ are compact, the horizontal maps factor through some finite stage~$X_\alpha$, and~$Y_\beta$.
We can take $\alpha=\beta$, increasing them if necessary.
By further increasing~$\alpha$ we can make the above diagram commutative for $s = \alpha$.
Since $X_\alpha\to Y_\alpha$ is an acyclic fibration, we have a lifting $B\to X_\alpha$, which gives a lifting of the original diagram
after postcomposing with $X_\alpha\to X_\infty$.
\ppar
\refit{we.chain.colimits}: We may assume that $\C$ is such that its generating cofibrations have compact (co)domains.
Suppose $Qf\to f$ is a cofibrant replacement of~$f$ in the projective structure.
Part~\refit{cof.sequential} shows that the transfinite composition of~$Qf$ is a weak equivalences, whereas part~\refit{af.sequential} shows that the filtered
colimit of the maps $QX_i\to X_i$ and $QY_i\to Y_i$ is a weak equivalence.
\refit{we.transfinite.composition} is a particular case of~\refit{we.chain.colimits}.
\xpf

The notion of h-cofibrations due to Batanin and Berger recalled below is the basis of (symmetric) h-monoidality
(Definitions \ref{defi--i.monoidal}, \ref{defi--symmetric.i.monoidal}), which is a key condition in the admissibility results of a subsequent paper \cop{\refth{O.Alg}}.
The key point of h-cofibrations is that in left proper model categories h-cofibrations coincide with maps along which cobase changes are homotopy cobase changes.
In fact, many results of this paper still hold for nonproper model categories if we replace h-cofibrations with such maps,
however, our main supply of h-cofibrations comes from h-monoidal categories, which are automatically left proper.

\defi \label{defi--i.cofibration}
\label{defi--h.cofibration}
\cite[Definition 1.1]{BataninBerger:Homotopy}
A map $f\colon X \r X'$ in a model category~$\C$ is an \emph{h-cofibration} if for any pushout diagram
$$\xymatrix @-1pc{
X \ar[d]_f \ar[r] & A \ar[d] \ar[r]^g & B \ar[d] \\
X' \ar[r] & A' \ar[r]^{g'} & B'
}$$
with a weak equivalence~$g$, $g'$~is also a weak equivalence.
An \emph{acyclic h-cofibration} is a map that is both an h-cofibration and a weak equivalence.
\xdefi

\exam
\label{exam--sSet.h.cofibrations}
In the category $\sSet$, equipped with its standard model structure, a map is an (acyclic) cofibration if and only if it is an (acyclic) h-cofibration.
By \refle{i.cofibrations}\refit{h.cofibration.containment}, we only need to prove the if-part.
Suppose a noninjective map $f\colon A \to B$ is an h-cofibration.
Then $A$~has two nondegenerate simplices $a, a' \in A_n$ with $f(a) = f(a')$.
Since any cofibration is an h-cofibration and h-cofibrations are stable under composition by \ref{lemm--i.cofibrations}\refit{h.i.cofibrations.pushout}, we may first replace~$A$ by the union of all faces of~$a$ and~$a'$ and then by $S^n \vee S^n$,
using the pushout along the map $A \r S^n \vee S^n$ collapsing all proper faces of~$a$ and~$a'$ to the base point.
The pushout of $B \sqcup_{S^n \vee S^n} S^n$ (using the obvious collapsing map) is isomorphic to~$B$.
If $B$ was also the homotopy pushout, there was a homotopy fiber square of derived mapping spaces
$$\xymatrix @-1pc{
\hMap (S^n \vee S^n, K(\Z, n)) & &  \hMap (S^n, K(\Z, n)) \ar[ll] \\
\hMap (B, K(\Z, n)) \ar[u]^{f^{*}} & & \hMap (B, K(\Z, n)) \ar[u] \ar[ll]^{\id},
}$$
contradicting the fact that the path components of these spaces are $\Z \oplus \Z$, $\Z$, and $\H^n(B, \Z)$, respectively.
\xexam

Usually, h-cofibrations form a strictly larger class than cofibrations, though.
We don't know an effective criterion characterizing h-cofibrations.

In the following lemma, we write $\cofh$ and $\Acofh$ for the classes of h-cofibrations and acyclic h-cofibrations.
We denote the class of maps $f$ such that all pushouts along $f$ are homotopy pushouts by $\cofi$.
Similarly, $\Acofi := \cofi \cap \we$.

\lemm
\label{lemm--i.cofibrations}
Suppose $\C$ is a model category.

\begin{enumerate}[(i)]
\item
\label{item--h.cofibration.containment}
We have inclusions
$$\cofi \subset \cofh, \qquad \AC \subset \Acofi \subset \Acofh.$$
Moreover, $\C$ is left proper if and only if $\cof \subset \cofi$ (equivalently, $\cof\subset\cofh$), in which case we have
$$\cofi = \cofh, \qquad \AC \subset \Acofi = \Acofh.$$

\item
\label{item--h.i.cofibrations.pushout}
The classes of $\cofi$, $\cofh$, and $\Acofi$ are stable under composition, pushouts, and retracts.
The same is true for $\Acofh$ if $\C$ is left proper.

\item
\label{item--acyclic.i.cofibrations}
The class $\Acofi$ consists precisely of the couniversal weak equivalences, i.e., weak equivalences that are stable under arbitrary pushouts.

\item \label{item--i.cofibrations.weakly.sat}
If weak equivalences are stable under colimits of chains (e.g., if $\C$ is pretty small, see \refle{sequential}\refit{we.chain.colimits}), then so are $\cofh$, $\cofi$, $\Acofi$, and $\Acofh$.
Thus, the first three among these classes are weakly saturated, and the same is true for $\Acofh$ if $\C$ is in addition left proper.
\end{enumerate}
\xlemm

\pf
We use the following well-known fact: given a weak equivalence $f\colon A \r B$, a pushout $f'\colon A' \r A' \sqcup_A B$ is a homotopy pushout if and only if $f'$ is a weak equivalence.
This follows from applying the 2-out-of-3-property of weak equivalences to
$$f'\colon A' \stackrel \sim \r A' \sqcup_A^h A \stackrel \sim \r A' \sqcup_A^h B \r A' \sqcup_A B.$$
(The first map is always a weak equivalence by computing $A' \sqcup^h_A A$ as $QA' \sqcup_{QA} QA$,
where $QA \r QA'$ is a cofibrant replacement of the map $A \r A'$, i.e., a cofibration with cofibrant source.)

Parts~\refit{h.cofibration.containment} and~\refit{h.i.cofibrations.pushout} are due to Batanin and Berger \cite[Proposition~1.6, Lemmas 1.2, 1.3,~1.7]{BataninBerger:Homotopy}.
(The characterization of left properness in~\refit{h.cofibration.containment} is formulated for $\cofh$ instead of $\cofi$ in loc.~cit., but holds as stated above by the above fact.)
\refit{h.i.cofibrations.pushout} also implies the ``$\Rightarrow$'' implication of~\refit{acyclic.i.cofibrations}.
The converse follows from the above fact.

\refit{i.cofibrations.weakly.sat}:
It is enough to treat the non-acyclic classes.
We first show the claim for $\cofh$, using the notation of \refle{sequential}.
For an object~$S$ under $X_\infty$, there is a functorial isomorphism $S \sqcup_{X_\infty} Y_{\infty} = \colim S \sqcup_{X_i} Y_i$.
Therefore, the pushout of a weak equivalence $s\colon S\to S'$ under $X_\infty$ along $f_\infty$ is the filtered colimit of the pushouts of $s \sqcup_{X_i} Y_i$.
Each of those is a weak equivalence since $f_i$ is an h-cofibration.
By assumption, their colimit is also a weak equivalence, so $f_\infty$ is an h-cofibration.
A similar argument works for $\cofi$ by commuting filtered homotopy colimits and homotopy pushouts.
\xpf

\lemm
\label{lemm--i.cofibrations.detect}
Suppose $G\colon \D \r \C$ is a functor between model categories that creates weak equivalences (for example, if the model structure on $\D$ is transferred from $\C$).
If $G$ preserves pushouts along a map $d \in \Mor(\D)$ and $G(d)$ is an (acyclic) h-cofibration, then $d$~is an (acyclic) h-cofibration.
\xlemm

\pf
Given a pushout $f'$ in $\D$ of a weak equivalence~$f$ under $\dom(d)$, we apply~$G$ and get a pushout in~$\C$.
As $G(d)$ is an h-cofibration, $G(f')$ is a weak equivalence, hence $f'$ is a weak equivalence and therefore $d$~is an h-cofibration.
The acyclic part is similar, using that $G$~detects weak equivalences.
\xpf

\numberwithin{equation}{subsection}

\section{Monoidal model categories}

In this section, we study certain properties of monoidal model categories.
We begin with a discussion of the pushout product, which is the relevant monoidal structure on arrows in a monoidal category.
We then recall the concepts of h-monoidality (due to Batanin and Berger) and flatness (due to Hovey) and establish a weak saturation property.
In the case of a symmetric monoidal model category, these notions will be refined in \refchap{symmetric}.

\subsection{The pushout product}
\label{sect--arrows}
\label{sect--pushout.product}

In this section, we define an endofunctor~$\Ar$ on the bicategory of cocomplete monoidal categories, cocontinuous strong monoidal functors, and monoidal natural transformations.
Roughly speaking, $\Ar$ sends a category~$\C$ to its category of morphisms equipped with a new monoidal structure, the {\it pushout product}.
The underlying category of~$\Ar(\C)$ is the category of functors~$\Fun(2,\C)$, where $2:=\{0\to1\}$ is the arrow category.
Its objects are morphisms in~$\C$ and its morphisms are commutative squares in~$\C$.
If $\C$ is (co)complete, then $\Ar(\C)$ is also (co)complete, because (co)limits in categories of functors are computed componentwise.
In this section we study the monoidal structure of $\Ar(\C)$ given by the pushout product and the projective model structure on $\Ar(\C)$.

\defi \label{defi--Ar}
Given a cocomplete monoidal category~$\C$,
its (cocomplete) category~$\Ar(\C)$ of morphisms can be endowed with a monoidal structure (the {\it pushout product\/}) as follows.
Interpret an object in~$\Ar(\C)$ as a functor $2\to \C$.
A finite family~$f\colon I\to \Ar(\C)$ of objects in~$\Ar(\C)$ (i.e., morphisms $f_i\colon X_i\to Y_i$ in~$\C$)
gives a functor~$2^I\to \C^I\to \C$, where $\C^I\to \C$ is the monoidal product on~$\C$.
We interpret this functor as a cocone on the category $2^I\setminus\{1^I\}$ (observe that $1^I$ is the terminal object of the category~$2^I$)
and the monoidal product of~$f$ is defined to be the universal map $\bigpp f_i\colon\bigppdom f_i\to\bigotimes_i Y_i$ associated to this cocone, interpreted as an object in~$\Ar(\C)$.
This defines a monoidal structure on $\Ar(\C)$.
\ppar
For example, the pushout product of two morphisms $f_1$~and~$f_2$ is
$$f_1\pp f_2\colon f_1 \ppdom f_2 = X_1 \t Y_2 \sqcup_{X_1 \t X_2} Y_1 \t X_2\to Y_1 \t Y_2.$$
We obtain a bifunctor
$$\pp\colon\Ar(\C)\x\Ar(\C)\to\Ar(\C).$$
\xdefi

\rema \label{rema--symmetric.closed}
If $(\C, \t)$ is braided or symmetric, then so is $(\Ar(\C), \pp)$.
Moreover, if $\t$ preserves colimits of a certain type (e.g., sifted colimits) in one or both variables, then so does $\pp$.
For example, if $\C$ is a closed monoidal category, then so is $\Ar(\C)$, with the internal hom $\Hom(f_1,f_2)$
(which one can call the {\it pullback hom\/} from~$f_1$ to~$f_2$) being the morphism
$\Hom(Y_1,X_2)\to\Hom(Y_1,Y_2)\times_{\Hom(X_1,Y_2)}\Hom(X_1,X_2)$.
For brevity of the exposition, we only spell out the nonsymmetric, nonclosed case in the sequel.
\xrema

\prop \label{prop--Ar.functorial}
If $F$ is a cocontinuous strong monoidal functor between cocomplete monoidal categories, then so is $\Ar(F)$.
\xprop

\pf
The functor~$\Ar(F)$ is cocontinuous because colimits of diagrams are computed componentwise.
To prove strong monoidality, suppose $f\colon I\to\Ar(\C)$ is a finite family of objects in~$\Ar(\C)$.
The diagram
$$\xymatrix @ -.3pc{
  2^I \ar[r]^f \ar[d]^\id & \C^I \ar[r]^\t \ar[d]^{F^I} & \C \ar[d]^F \\
  2^I \ar[r]^{F(f)} & \D^I \ar[r]^\t & \D
}$$
is commutative, meaning the left square is strictly commutative
and the right square is commutative up to the canonical natural isomorphism coming from the monoidal structure on the functor~$F$.
The pushout product $\bigpp f$ is the universal map associated to the cocone $2^I\To1f\C^I\To1\otimes\C$ with the apex $1^I\in2^I$,
and similarly for~$\bigpp\Ar(F)(f)$.
Since~$F$ is cocontinuous, it preserves universal maps associated to cocones.
Thus, the image of the universal morphism associated to the cocone $2^I\to \C^I\to \C$
is also the universal morphism associated to the cocone $2^I\to \C^I\to \C\to \D$.
The latter cocone is canonically isomorphic to the cocone $2^I\to \D^I\to \D$, which is the cocone defining~$\bigpp\Ar(F)(f)$.
\xpf

\defi
A morphism in the category~$\Ar(\C)$ for some monoidal category~$\C$
is a {\it pushout morphism\/} if the corresponding commutative square in~$\C$ is cocartesian.
\xdefi

\prop \label{prop--pushout.product.pushout.morphisms}
For any cocomplete closed monoidal category~$\C$ pushout morphisms in~$\Ar(\C)$ are closed under the pushout product.
\xprop

\pf
A pushout morphism can be presented as a functor~$2\times2\to \C$,
where the first~$2$ is responsible for the morphism direction in~$\Ar(\C)$ and the second~$2$ is responsible for the morphism direction in~$\C$.
Schematically, we denote this by the commutative diagram 
$$\xymatrix @ -1.2pc{
  00 \ar[rr] \ar[d] & & 10 \ar[d] \\
  01 \ar[rr] & & 11. 
}$$
A finite family of pushout morphisms~$f\colon I\to\Mor(\Ar(\C))$
gives a functor $(2\times2)^I\to \C^I$, which we compose with the monoidal product $\C^I\to \C$ to obtain a functor $F\colon(2\times2)^I\to \C$.
Consider now the category~$\DC$ of all full subcategories~$A$ of~$(2\times2)^I$ that are {\it downward closed\/} (or {\it convex\/} in the sense of Goodwillie):
if $Y\in A$ and $X\to Y$ is a morphism in~$(2\times2)^I$, then also $X\in A$.
Morphisms in~$\DC$ are inclusions of subcategories.
Taking the colimit of the functor~$F$ restricted to the given full subcategory~$A$ yields a cocontinuous functor $Q\colon\DC\to \C$.
In particular, the set of all inclusions~$A\to B$ in~$\DC$ that are mapped to isomorphisms by~$Q$ forms a subcategory of~$\DC$ closed under cobase changes of the underlying sets.
\ppar
Suppose that $B\in\DC$ is obtained from~$A\in\DC$ by adding an element $W\times11$ and taking the downward closure,
where $W\in(2\times2)^{I\setminus i}$ for some $i\in I$ is such that $W\times\{00,01,10\}\subset A$.
The resulting inclusion $A\to B$ gives an isomorphism after we apply~$Q$ because the commutative square~$2\times2\To2{\times W}(2\times2)^I\To2F\C$
is a cocartesian square because each $f_i$ is a cocartesian square and the monoidal product with a fixed object preserves cocartesian squares.
This uses the closedness of the monoidal product.
\ppar
Consider the commutative square in $\DC$,
whose right entries are obtained by taking the left entries, replacing~0 in the first components by~1, and downward closing:
$$\xymatrix @-1.2pc{
  \{00,01\}^I\setminus\{01\}^I \ar[rr] \ar[d] & & \{00,01,10,11\}^I\setminus\{01,11\}^I \ar[d] \\
\{00,01\}^I \ar[rr] & &  \{00,01,10,11\}^I.}$$
The pushout product $\pp f_i$ is obtained by applying~$Q$ to the following map:
$$\{00,01,10,11\}^I\setminus\{01,11\}^I \sqcup_{\{00,01\}^I\setminus\{01\}^I} \{00,01\}^I\to\{00,01,10,11\}^I.$$
We present this morphism in~$\DC$ as a composition of pushouts of generating maps explained in the previous paragraph,
which implies that the map itself is sent to an isomorphism by~$Q$.
Such a presentation can be obtained by using the rule explained above to add all elements of~$\{01,11\}^I\setminus\{01\}^I$ to the source by induction on the number of~$11$'s.
If there are no~$11$'s, the element $\{01\}^I$ belongs to the bottom left corner, proving our claim.
By induction, assuming that all tuples with less than~$k$ elements equal to~$11$ have already been added,
take any tuple with exactly $k$~components equal to~$11$ and observe that by replacing this component with $00$,~$01$, or~$10$ we obtain a tuple
already present in our set.
Thus, we can also add the tuple under consideration to our set.
\xpf

The elementary proof of the following lemma is left to the reader.
Together with \refpr{pushout.product.pushout.morphisms}, it can be rephrased by saying that $x \pp -$ preserves finite cellular maps.

\lemm
\label{lemm--pp.compositions}
Given two composable maps $y$~and~$z$, and another map~$x$, $x \pp (y \circ z)$ is the composition of the pushout of $x \pp z$ along $x\ppdom z\to x\ppdom(y\circ z)$, followed by $x \pp\nobreak y$.
\xlemm

We now extend the formation of arrow categories to monoidal model categories.
A {\it strong monoidal left Quillen functor\/} between monoidal model categories
is a left Quillen functor~$F$ that is also equipped with the structure of a strong monoidal functor, i.e., functorial isomorphisms $F(X \t Y) \cong F(X) \t F(Y)$ compatible with the unit and
associativity of $\t$.
Monoidal model categories, strong monoidal left Quillen functors, and monoidal natural transformations form a bicategory.
(As in \refre{symmetric.closed}, there are obvious variants for (symmetric) monoidal model categories, which we will not spell out explicitly.)

The following proposition was shown independently by Hovey under the additional assumption that $\C$ is cofibrantly generated \cite[Proposition~3.1]{Hovey:Smith}.

\prop \label{prop--Ar.model}
The functor~$\Ar$ (\refde{Ar}) descends to the bicategory of monoidal model categories, as described in the proof below.
\xprop

\pf
Given a closed monoidal model category~$\C$, the monoidal category~$\Ar(\C)$ is complete and cocomplete.
We equip $\Ar(\C)$ with the projective model structure, which coincides with the Reedy model structure, where the nonidentity arrow~$0\to1$ in~$2$ is declared to be positive.
In particular, the projective model structure on~$\Ar(\C)$ exists.
Fibrations and weak equivalences are defined componentwise.
(Acyclic) cofibrations $f\colon g\to h$ are commutative squares
$$\cd{W&\mapright p&Y\cr\mapdown g&&\mapdown h\cr X&\mapright q&Z\cr}$$
such that $p$ and the universal map $Y\sqcup_W X\to Z$ are both (acyclic) cofibrations, hence $q$ is also an (acyclic) cofibration.
In particular, cofibrant objects in~$\Ar(\C)$ are morphisms~$g\colon W\to X$ such that $W$ is cofibrant and $g$ is a cofibration in~$\C$.
\ppar
We now prove the pushout product axiom for $\Ar(\C)$ from the one of $\C$ (\refde{monoidal.model.category}).
Actually, we show that the pushout product of a finite nonempty family~$f\colon I\to\Mor(\Ar(\C))$ of cofibrations
in~$\Ar(\C)$ is a cofibration,
and if one of the cofibrations is acyclic, then the resulting cofibration is also acyclic.
The infrastructure of the following proof is the same as in the proof of \refpr{pushout.product.pushout.morphisms}.
Just like there we get a functor $F\colon(2\times2)^I\to \C$ and a cocontinuous functor $Q\colon\DC\to \C$.
Let
$$\cd{A&\mapright {}&A'\cr\mapdown a&&\mapdown {a'}\cr B&\mapright{} &B'\cr}$$
be a cocartesian square in $\DC$, i.e., $B' = A' \cup_A B$.
If $Q(a)$ is a cofibration, then so is $Q(a')$.
Suppose that for every~$i\in I$ we select one of the morphisms $\{00\}\to\{00,10\}$ or $\{00,01,10\}\to\{00,01,10,11\}$ in~$\DC(2\times 2)$.
Then the pushout product of these morphisms belongs to the above subcategory because of the pushout product axiom for~$\C$.
The first morphism above expresses the fact that the top arrow of a cofibration in~$\Ar(\C)$ is itself a cofibration
and the second morphism corresponds to the canonical map from the pushout to the bottom right corner, which is also a cofibration.
The pushout product mentioned above always has the form $A\setminus\{x\}\to A$, where the individual components of~$x$ are $10$ respectively~$11$,
according to the choice made above.
\ppar
The pushout product of~$f$ is the functor~$Q$ applied to the commutative square
$$\cd{\{00,01,10,11\}^I\setminus\{10,11\}^I\setminus\{01,11\}^I&\to&\{00,01,10,11\}^I\setminus\{01,11\}^I\cr
\downarrow&&\downarrow\cr
\{00,01,10,11\}^I\setminus\{10,11\}^I&\to&\{00,01,10,11\}^I.\cr}$$
It remains to prove that $Q$ applied to the top map and the map from the pushout of the left and top arrows (i.e., the union of all corners except for the bottom right corner)
to the bottom right corner is a cofibration.
We present the morphism in~$\DC$ under consideration as a composition of pushouts of generating maps explained in the previous paragraph.
This implies that the map itself is sent by $Q$ to a cofibration.
\ppar
For the top map, such a presentation can be obtained by using the rule explained above to add all elements of~$\{10,11\}^I\setminus\{11\}^I$
to the source by induction on the number of~$11$'s.
Assume that all tuples with less than~$k$ $11$'s have already been added and take any tuple with exactly $k$ $11$'s.
By applying the rule explained in the previous paragraph
to the family of maps that are either $\{00\}\to\{00,10\}$ if the corresponding component is~$10$
or $\{00,01,10\}\to\{00,01,10,11\}$ if the corresponding component is~$11$
we can conclude that the tuple under consideration can be added to our set.
\ppar
For the map from the pushout of the top and left arrows to the bottom right corner observe that we only need to add the element~$\{11\}^I$,
which is possible because the conditions for the corresponding rule are satisfied.
\ppar
For acyclic cofibrations observe that the rule in the previous paragraph now guarantees that the resulting map is an acyclic cofibration after we apply~$Q$,
precisely because the pushout product in~$\C$ of a family of cofibrations, at least one of which is acyclic, is again an acyclic cofibration.
The rest of the proof is exactly the same, because the category of acyclic cofibrations is also closed under pushouts.
\ppar
Finally,~$\Ar$ descends to strong monoidal left Quillen functors:
if $F\colon \C\to \D$ is such a functor, then the induced functor $\Ar(F)\colon\Ar(\C)\to\Ar(\D)$ is cocontinuous and strong monoidal (\refpr{Ar.functorial}).
It is a left Quillen functor because $F$ preserves (acyclic) cofibrations and pushouts.
\xpf

\subsection{h-monoidality and flatness}
\label{sect--i.monoidal.flat}

\defi \label{defi--monoidal.model.category}
A \emph{(symmetric) monoidal model category} is a closed (symmetric) monoidal category~$\C$ such that $\t \colon \C \x \C \r \C$
is a left Quillen bifunctor \cite[Definition~4.2.6]{Hovey:Model}.
This is also referred to as the \emph{pushout product axiom}.
(Thus, unlike \cite[Definition 4.2.6]{Hovey:Model}, we do not require the unit axiom,
which asks that $(\mcr (1)\to 1)\t X$ is a weak equivalence for any cofibrant object~$X$,
because it is a special case of flatness, see \refde{flat}).

We say $\C$ satisfies the {\it monoid axiom\/} if the class $\ws{\C \t \AC_\C}$ consists of weak equivalences in $\C$ \cite[Definition~3.3]{SchwedeShipley:Algebras}.
\xdefi

In this section, we discuss the notion of h-monoidality and flatness of a monoidal model category $\C$.
The notion of h-monoidality was introduced by Batanin and Berger \cite[Definition 1.7]{BataninBerger:Homotopy}.
Essentially, h-monoidality ensures that category of modules over some monoid $R \in \C$ carries a model structure.
This statement is referred to as the \emph{admissibility} of the monoid $R$.
The admissibility of monoids is also guaranteed by the monoid axiom \cite[Theorem 4.1]{SchwedeShipley:Algebras}, which is a combination of two weak saturation properties, namely weak saturation by transfinite compositions and by pushouts.
In this paper, we focus on admissibility conditions using pretty smallness and h-monoidality, which individually govern the homotopical behavior of transfinite compositions and of (certain) pushouts, respectively.
The standard basic model categories are h-monoidal (\refchap{examples}).
On the other hand, h-monoidality is very robust since is stable under transfer and localization (\ref{prop--transfer.monoidal}\refit{transfer.i.monoidal}, \ref{prop--monoidal.localization}\refit{localization.i.monoidal}).
We don't know a similar statement for the monoid axiom (without the detour via pretty smallness and h-monoidality).

\defi \label{defi--i.monoidal}
A class $S$ of (acyclic) cofibrations in a monoidal category $\C$ is \emph{(acyclic) h-monoidal} if for any object~$C\in\C$ and any $s \colon S_1 \r S_2$ in $S$, the map
$$C \t s\colon C \t S_1 \r C \t S_2$$
is an (acyclic) h-cofibration (\refde{i.cofibration}).
The category $\C$ is \emph{h-monoidal} if the classes of (acyclic) cofibrations are (acyclic) h-monoidal.
\xdefi

Any h-monoidal model category is left proper \cite[Lemma 1.8]{BataninBerger:Homotopy}.
The following immediate consequence of \refle{i.cofibrations}\refit{acyclic.i.cofibrations},~\refit{i.cofibrations.weakly.sat}
is very similar to \cite[Proposition~2.5]{BataninBerger:Homotopy}.

\lemm
\label{lemm--i.monoidal.monoid.axiom}
If $\C$ is a h-monoidal model category whose acyclic h-cofibrations are stable under transfinite compositions (for example, $\C$ is pretty small),
then $\C$ satisfies the monoid axiom.
\xlemm

We now define flatness, which is the main condition in rectification of modules over monoids.
Its symmetric strengthening, symmetric flatness, plays the corresponding role for algebras over symmetric operads \cop{\refth{rect.colored.operad}}.

\defi \label{defi--flat}
A class $S$ of cofibrations in a monoidal model category $\C$ is \emph{flat} if for all weak equivalences $y \colon Y_1 \r Y_2$ in $\C$ and all $s \colon S_1 \r S_2$ in $S$, the following map is a weak equivalence:
$$y \pp s \colon Y_2 \t S_1 \sqcup_{Y_1 \t S_1} Y_1 \t S_2 \r Y_2 \t S_2.$$
The category $\C$ is \emph{flat} if the class of all cofibrations is flat.
\xdefi

For example, if $\C$ is flat then for any cofibrant object $X \in \C$ and any weak equivalence $y \in \C$, the map $y \t X$ is a weak equivalence.
In this slightly weaker form, flatness is independently due to Hovey \cite[Definition~2.4]{Hovey:Smith}.
Actually, the notion appears already in \cite[Theorem~4.3]{SchwedeShipley:Algebras}.
We use the above slightly stronger definition since it is stable under weak saturation of $S$ (\refth{power.weakly.saturated}\refit{flat.weakly.sat}).
This is useful to show the stability of flatness under transfer (\refpr{transfer.monoidal}\refit{transfer.flat}) and localization (\refpr{monoidal.localization}\refit{localization.flat}).

If all objects in a monoidal model category $\C$ are cofibrant then $\C$ is left proper \cite[Corollary~13.1.3]{Hirschhorn:Model}, quasi-tractable, h-monoidal \cite[Lemma 1.8]{BataninBerger:Homotopy}, and flat \cite[Remark~3.4]{SchwedeShipley:Algebras}.
In general, though, we avoid cofibrancy hypotheses where possible, in particular, we do not in general assume that the monoidal unit~$1$ is cofibrant.

We finish this section with two weak saturation properties.
A slightly weaker statement than \refth{power.weakly.saturated}\refit{flat.weakly.sat} is independently due to Hovey \cite[Theorem~A.2]{Hovey:Smith}.
The following lemma is the basis of the interaction of h-monoidality and flatness, see, for example, the proof of \ref{theo--power.weakly.saturated}\refit{flat.weakly.sat}.

\lemm \label{lemm--preparation.flat.weakly.sat}
Let $\C$ be a left proper monoidal model category.
Let $$\xymatrix{
A \ar[r] \ar[d]^a & B \ar[d]^b \\
A' \ar[r] & B'
}$$
be a cocartesian square in $\C$.
Let $y \colon Y \r Y' \in \C$ be any morphism such that $y \pp a$ is a weak equivalence in~$\C$, and both $Y \t a$ and $Y' \t a$ are h-cofibrations (\refde{h.cofibration}).
Then $y \pp b$ is a weak equivalence.
\xlemm
\pf
Consider the commutative diagram
$$\xymatrix{
Y \t A \ar[r]^{y \t A} \ar[d]_{Y \t a}  &
Y' \t A \ar[d]^\alpha \ar[drr]^{Y' \t a}  \\
Y \t A' \ar[r] &
y \ppdom a \ar[rr]_(.6){y \pp a} & &
Y' \t A'.
}$$
As usual, $\ppdom$ denotes the domain of the pushout product $\pp$.
By assumption, $Y \t a$ is an h-cofibration, hence so is $\alpha$ by \refle{i.cofibrations}\refit{h.i.cofibrations.pushout}.
Likewise, $Y' \t a$ is an h-cofibration.
Hence, the top square and outer rectangle in the diagram below are homotopy pushouts (\refle{i.cofibrations}\refit{h.cofibration.containment})
and so is the bottom square.
The map $y \pp b$ is therefore also a weak equivalence:
$$\xymatrix{
Y' \t A \ar[r] \ar[d]_\alpha^{\text{h-cofib.}} \ar@/_2pc/[dd]_{\text{h-cofib.}}&
Y' \t B \ar[d]\\
y \ppdom a \ar[r] \ar[d]^{y \pp a}_\sim &
y \ppdom b \ar[d]_{y \pp b}
\\
Y' \t A' \ar[r] &
Y' \t B'.
}$$
\xpf

\theo \label{theo--power.weakly.saturated}
Let $\C$ be a left proper, pretty small monoidal model category.
We say some property of a class $S$ of morphisms in $\C$ is \emph{stable under saturation} if it also holds for the weak saturation $\ws S$.
\begin{enumerate}[(i)]
\item \label{item--i.monoidal.weakly.sat}
The property of $S$ of being (acyclic) h-monoidal is stable under saturation.

\item \label{item--flat.weakly.sat}
If $S$ is h-monoidal, then the flatness of $S$ is stable under saturation.
In particular, if some class of generating cofibrations in $\C$ is flat and h-monoidal, then $\C$ is flat.
\end{enumerate}
\xtheo

\rema
\label{rema--h.monoidal.weakly.sat}
The proof below shows that the left properness assumption on $\C$ is not necessary to show the non-acyclic part in~\refit{i.monoidal.weakly.sat}.
This also implies that $\C$ is h-monoidal (asssuming pretty smallness, but not left properness) if some sets $I_\C$ (respectively, $J_\C$) of generating (acyclic) cofibrations are (acyclic) h-monoidal, since any h-monoidal model category is left proper.
\xrema

\pf
\refit{i.monoidal.weakly.sat}:
This follows from
\refle{i.cofibrations}\refit{i.cofibrations.weakly.sat} and the cocontinuity of $\t$.

\refit{flat.weakly.sat}:
For a weak equivalence $y\colon Y \r Y'$ in $\C$ and any $s \in S$, $y \pp s$ is a weak equivalence by assumption.
By h-monoidality of $S$, $Y \t s$ and $Y' \t s$ are h-cofibrations.
Thus, for any pushout $s'$ of $s$, $y \pp s'$ is a weak equivalence by \refle{preparation.flat.weakly.sat}.
For a transfinite composition $s_\infty$ of maps $s_i$, $y \pp s_\infty$ is the transfinite composition of $y \pp s_i$ by preservation of filtered colimits in the second variable.
Therefore it is again a weak equivalence using pretty smallness (\refle{sequential}).
As usual, retracts are clear.
\xpf

\coro
\label{coro--acyclic.h.monoidal}
Any quasi-tractable, flat, pretty small monoidal model category $\C$ that satisfies the nonacyclic part of h-monoidality also satisfies the acyclic part of h-monoidality.
\xcoro
\ppar

\pf
The nonacyclic part of h-monoidality implies left properness by \refle{i.cofibrations}\refit{h.cofibration.containment}.
We can thus apply \refth{power.weakly.saturated}\refit{i.monoidal.weakly.sat}: it is enough to show that for any object $X$ and any generating acyclic cofibration $s\colon S \r S'$,
the map $X \t s$ is a weak equivalence.
By quasi-tractability we may assume that $S$ (and $S'$) is cofibrant.
Pick a cofibrant replacement $q\colon \mcr X \r X$.
Then $X \t s$ is a weak equivalence since $\mcr X \t s$ is one
(by the pushout product axiom) and $q \t S$ and $q \t S'$ are weak equivalences in~$\C$ by flatness.
\xpf

\rema
\label{rema--monoidality.acyclic.part}
One shows similarly that any quasi-tractable flat model category that satisfies the nonacyclic part of the pushout product axiom also satisfies the acyclic part of the pushout product axiom.
\xrema


\section{Symmetricity properties}
\label{chap--symmetric}

In this section we study three properties of a symmetric monoidal model category $\C$: symmetric h-monoidality, symmetroidality and symmetric flatness.
As the name indicates, these involve the formation of pushout powers, i.e., expressions of the form
$$\bigpp_n f := f^{\pp n} := \underbrace{f \pp \cdots \pp f}_{\hbox{$n$~times}}.$$
After settling preliminaries on objects with a finite group action, these properties are defined in~\refsect{definitions}.
The main result of~\refsect{symmetric.weak.saturation} is \refth{symmetric.weakly.saturated} which shows the stability of these notions under weak saturation.
This is a key step in showing that the properties also interact well with transfer and localization of model structures.
Examples of model categories satisfying these properties are given in \refchap{examples.symmetric}.

\subsection{Objects with a finite group action}
\label{sect--finite.group.action}

We first examine model-theoretic properties of objects with an action of a finite group, for example, the permutation action of $\Sigma_n$ on $f^{\pp n}$.
Given a finite group~$G$, considered as a category with one object, and any category $\C$, define
$$G\C:=\Fun(G,\C).$$
This is the category of objects in~$\C$ with a left $G$-action.
It is symmetric monoidal if $\C$ is, by letting $G$ act diagonally on the monoidal product.
Given some $X\in G\C$ and any subgroup $H\subset G$, we write $X_H=\colim_H X$ for the coinvariants.

For any $X \in \C$ we define $G/H\cdot X := \coprod_{G/H}X \in G\C$ on which $G$-acts by the left $G$-action on~$G/H$.
More generally, given any $X\in H\C$ and any morphism of groups $H \r G$, we define
$G\cdot_H X:=(G\cdot X)_H$, where $H$ acts on the right on~$G$ and on the left on~$X$.

\lemm \label{lemm--dirty}
Suppose $\C$~is a cocomplete category and $H$ is a subgroup of a finite group $G$.
Any choice of a partition $G = \coprod_i H \cdot g_i$ of $G$ into $H$-cosets
induces a natural isomorphism
$$\varphi (G \cdot_H -)\to (G/H) \cdot \varphi (-)$$
of functors $H\C\to\C$, where $\varphi$ denotes the forgetful functor to~$\C$.
\xlemm

\pf
The canonical projection $G \cdot \varphi X\to G/H \cdot \varphi X$ factors over $\varphi(G \cdot_H X)$.
Conversely, given $g \in G$, the partition gives a unique $h \in H$ and $i$ such that $g = h g_i$.
Define $G/H \cdot \varphi X\to G \cdot_H \varphi X$ by $x_{gH} \mapsto (h^{-1} x)_{g_i}$.
\xpf

\prop \label{prop--projective.G}
Suppose $\C$ is a cofibrantly generated model category.
The category $G\C$ carries the \emph{projective model structure}, denoted $G^\proj\C$,
whose weak equivalences and fibrations are precisely those maps in~$G\C$ that are mapped to weak equivalences, respectively, fibrations in~$\C$ by the forgetful functor~$G\C\to\C$.
The cofibrations of~$G^\proj\C$ are generated by the maps of the form $G\cdot f$, where $f$ runs over generating cofibrations of~$\C$.
\ppar
Given a morphism of groups $H\to G$, there is a Quillen adjunction
$$G\cdot_H- : H^\proj\C \leftrightarrows G^\proj\C : R,\eqlabel{Quillen.adjunction.H.G}$$
where the right adjoint functor is the restriction.
\ppar
Finally, suppose $\C$ is a symmetric monoidal model category.
Given two groups $G$~and~$H$, the monoidal product on~$\C$ induces a left Quillen bifunctor
$$G^\proj\C \times H^\proj\C \to(G\times H)^\proj\C.\eqlabel{Quillen.bifunctor.projective}$$
\xprop

\pf
The existence of this model structure is standard, see, for example, Hirschhorn \cite[Theorem~11.6.1]{Hirschhorn:Model}.
The adjunction~\refeq{Quillen.adjunction.H.G} is seen to be a Quillen adjunction by looking at the right adjoint.
The functor~\refeq{Quillen.bifunctor.projective} is a left Quillen bifunctor because $(G\cdot I_\C)\pp(H\cdot I_\C)=(G\times H)\cdot(I_\C\pp I_\C)\subset(G\times H)\cdot\cof_\C$,
using the cocontinuity and monoidality of the functor~$G\cdot-$ and the pushout product axiom for~$\C$.
\xpf

\subsection{Definitions of symmetricity properties}
\label{sect--definitions}

We now define three properties of (morphisms in) a symmetric monoidal model category~$\C$: symmetric flatness, symmetric h-monoidality and symmetroidality.
They are appropriate strengthenings of flatness (\refde{flat}), h-monoidality (\refde{i.monoidal}) and the pushout product axiom.
Symmetric flatness is the key condition required to obtain a rectification result for operadic algebras \cop{\refth{cow.Coll}}.
Approximately, it says that for any cofibrant object $X \in \C$, the map
$$y \t_{\Sigma_n} X^{\t n}\colon Y \t_{\Sigma_n} X^{\t n} \r Y' \t_{\Sigma_n} X^{\t n}$$
is a weak equivalence for any weak equivalence $y\colon Y \r Y'$ in $\Sigma_n \C$.
Slightly more accurately, the definition is phrased in terms of more general cofibrations $s$ using instead
$$y \pp_{\Sigma_n} s^{\pp n}.$$
For $s \colon \emptyset \r X$ this gives back the previous expressions.
In order to ensure that the three symmetricity properties are stable under weak saturation (\refth{symmetric.weakly.saturated}), we actually define them for a class of morphisms instead of a single morphism.
In such cases, we use the following notational conventions.

\defi \label{defi--multi}
Let $v := (v_1, \ldots, v_e)$ be a finite family of morphisms.
For any sequence of nonnegative integers $n := (n_i)_{i \leq e}$, we write $\Sigma_n := \prod \Sigma_{n_i}$, $v^{\pp n} := v_1^{\pp n_1}  \pp \cdots \pp v_e^{\pp n_e}$,
and $v^{\otimes n}:=v_1^{\otimes n_1}\otimes\cdots\otimes v_e^{\otimes n_e}$.
We write $m\le n$ if $m_i\le n_i$ for all~$i$
and $m<n$ if $m\le n$ and $m\ne n$.
Given a class $S$ of morphisms, we write $v \subset S$ if all $v_i$ are in $S$.
Given another sequence of integers $(m_i)_{i=1}^e$, we write $mn:=\sum m_i n_i$ and $\Sigma_m^n := \prod \Sigma_{m_i}^{n_i}$ and $\Sigma_n \rtimes \Sigma_m^n := \prod \Sigma_{n_i} \rtimes \Sigma_{m_i}^{n_i}$.
\xdefi

\defi \label{defi--eq.flat} \label{defi--symmetric.flat}
A class $S$ of cofibrations in $\C$ is called \emph{symmetric flat with respect to some class $\Y = (\Y_n)$} of morphisms $\Y_n \subset \Mor \Sigma_n \C$ if
$$y \pp_{\Sigma_n} s^{\pp n} := (y \pp s^{\pp n})_{\Sigma_n}$$
is a weak equivalence in $\C$ for any $y \in \Y_n$, any finite multi-index $n \ge 1$ and any $s \in S$.
We say $S$ is \emph{symmetric flat} if it is symmetric flat with respect to the classes $\Y_n = (\we_{\Sigma_n^\proj \C})$ of projective weak equivalences (i.e., those maps in $\Sigma_n \C$ which are weak equivalences after forgetting the $\Sigma_n$-action).
We say $\C$ is \emph{symmetric flat} if the class of cofibrations is symmetric flat.
\xdefi

\exam
\label{exam--no.multiindices}
A class $S$ is symmetric flat (i.e., with respect to $\we_{\Sigma_n^\proj \C}$) if and only if $y \pp_{\Sigma_n} s^{\pp n}$ is a weak equivalence for a single map $s \in S$, i.e., no multi-indices are necessary in this case.
The reader is encouraged to mainly think of this case.
\xexam

The notions of symmetric h-monoidal maps (respectively, symmetroidal maps) presented next are designed to ultimately address the (strong) admissibility of operads (\cop{\refth{O.Alg}}).

\defi
\label{defi--symmetric.i.monoidal}
A class $S$ of morphisms in a symmetric monoidal category $\C$ is called {\it(acyclic) symmetric h-monoidal} if for any finite family~$s \subset S$ and any multi-index $n\ne0$,
and any object~$Y \in \Sigma_n \C$
the morphism $Y\otimes_{\Sigma_n} s^{\pp n}$ is an (acyclic) h-cofibration.
We say $\C$ is \emph{symmetric h-monoidal} if the class of (acyclic) cofibrations is (acyclic) symmetric h-monoidal.
\xdefi

The notion of power cofibrations presented next is due to Lurie \cite[Definition~4.5.4.2]{Lurie:HA} and Gorchinskiy and Guletski\u\i~\cite[Section~3]{GorchinskiyGuletskii:Symmetric}, who also introduced symmetrizable maps.

\defi
\label{defi--symmetroidal}
\label{defi--power.cofibration}
Let $\Y=(\Y_n)_{n > 0}$ be a collection of classes $\Y_n$ of morphisms in $\Sigma_n\C$, where $n > 0$ is any finite multi-index.
We suppose that for $y \in \Y_n$, $y \pp -$ preserves \emph{injective} (acyclic) cofibrations in $\Sigma_n\C$, i.e., those maps which are (acyclic) cofibrations in $\C$.

A class $S$ of morphisms in a symmetric monoidal category $\C$ is called \emph{(acyclic) $\Y$-symmetroidal} if for all multi-indices $n > 0$ and all maps $y \in \Y_n$, the morphism
$$y \pp_{\Sigma_n} s^{\pp n}$$
is an (acyclic) cofibration in~$\C$ for all $s \in S$.
If $\Y_n = \cof_{\Sigma_n^\inje \C}$, we say that $S$ is \emph{(acyclic) symmetroidal}.
For $\Y_n = \{ \emptyset \r 1_\C\}$, we say $S$ is \emph{(acyclic) symmetrizable}.

A map $f \in \C$ is called an {\it (acyclic) power cofibration\/} if the morphism $f^{\pp n}$ is an (acyclic) cofibration in~$\Sigma_n^\proj\C$ for all integers $n > 0$ (i.e., a projective cofibration with respect to the $\Sigma_n$-action).

The category $\C$ is called symmetric h-monoidal/$\Y$-symmetroidal/freely powered if the class of all (acyclic) cofibrations is (acyclic) symmetric h-cofibrant/(acyclic) $\Y$-symmetroidal/(acyclic) power cofibration.
\xdefi

\rema
\label{rema--symmetroidal}
\label{rema--power}
In the definition of power cofibrations, no multi-indices are necessary: for power cofibrations~$s_i$
and any multi-index $n = (n_i)$, $s^{\pp n} := \bigpp_i s^{\pp n_i}$ is a $\Sigma_n := \prod \Sigma_{n_i}$ projective cofibration by the pushout product axiom.

Unlike the definition of power cofibrations in \cite{Lurie:HA}, we exclude the case $n=0$, for this would require $1$ to be cofibrant, which is not always satisfied.
In fact, it is \emph{never} satisfied for the positive model structures on symmetric spectra which is a main motivating example for us \cite{PavlovScholbach:Spectra}.

We have the following implications (where symmetroidality is with respect to the classes $\Y_n$ of injective cofibrations in $\Sigma_n \C$):
$$\vcenter{\xymatrix{\text{power cofibration}  \ar@{=>}[r] &\text{symmetroidal map}  \ar@{=>}[r] \ar@{.>}[d]&\text{cofibration}  \ar@{=>}[d]\cr
&\text{symmetric h-cofibration} \ar@{=>}[r]&\text{h-cofibration}.\cr}}\eqlabel{implications}$$
The vertical implication holds if $\C$ is left proper.
The dotted arrow is not an implication in the strict sense unless all objects in~$\C$ are cofibrant.
A symmetroidal map $x$ is such that for all cofibrant objects $Y \in \Sigma_n^\inje \C$, the map $Y \t_{\Sigma_n} x^{\pp n}$ is a cofibration
and therefore (again if $\C$ is left proper) an h-cofibration.
Being a symmetric h-cofibration demands the latter for any object $Y\in\Sigma_n \C$.
Every power cofibration is a symmetrizable cofibration since the coinvariants $\Sigma_n^\proj\C\to\C$ are a left Quillen functor.
The implications in~\refeq{implications} are in general strict: in a monoidal model category $\C$ with cofibrant monoidal unit or, more generally, one satisfying the strong unit axiom, every object is h-cofibrant \cite[Proposition~1.17]{BataninBerger:Homotopy}, but of course not necessarily cofibrant.
In the category $\sSet$ of simplicial sets every cofibration is a symmetrizable cofibration, but not a power cofibration (see~\refsect{simplicial.sets}).

The homotopy orbit $\hocolim_{\Sigma_n} X^{\t n}$ can be computed by applying the derived functor of the either of the following two left Quillen bifunctors
to~$(1_\V,X^{\t n})$ \cite[Theorem~3.2 and Theorem~3.3]{Gambino:Weighted}:
$$\eqalignno{\Sigma_n^{\op,\inje}\V\x\Sigma_n^\proj\C & \To1\otimes \C,&\eqlabel{inj.proj}\cr
\Sigma_n^{\op,\proj}\V\x\Sigma_n^\inje\C & \To1\otimes \C.&\eqlabel{proj.inj}\cr}$$
Here $\V$ denotes the symmetric monoidal model category used for the enrichment and
the monoidal unit $1_\V \in \V$ is equipped with the trivial $\Sigma_n$-action.
If $\C$ is freely powered, then for any cofibrant object $X \in \C$, $X^{\t n}$ is projectively cofibrant, i.e., cofibrant in $\Sigma_n^\proj \C$.
Thus, the homotopy orbit is given by $(X^{\t n})_{\Sigma_n}$, provided that $1_\V$ is cofibrant \cite[Lemma~4.5.4.11]{Lurie:HA}.
However, most model categories appearing in practice are not freely powered, so that $X^{\t n}$ needs to be \emph{projectively} cofibrantly replaced to compute the homotopy colimit.
This is usually a difficult task.
On the other hand, when using~\refeq{proj.inj}, one needs to cofibrantly replace~$1$ in $\Sigma_n^{\op, \proj} \V$,
but no cofibrant replacement has to be applied to $X^{\t n}$, provided that $X$ is cofibrant in~$\C$.
This makes the second approach to computing homotopy colimits much more easily applicable.
This observation is used in \refle{preparatory} below, which in its turn is the key technical step in establishing the compatibility of symmetric h-monoidality and Bousfield localizations (\refth{symmetric.monoidal.localization}\refit{localization.symmetric.i.monoidal}).
\xrema

\subsection{Weak saturation of symmetricity properties}
\label{sect--symmetric.weak.saturation}

In this section, we provide a few elementary facts concerning the symmetricity notions defined in~\refsect{definitions}.
After this, we show the main theorem of this section (\ref{theo--symmetric.weakly.saturated}), which asserts that the symmetricity notions behave well with respect to weak saturation.

The following two results have a similar spirit: we show that symmetric flatness can be reduced to (projective) acyclic fibrations, and that the class $\Y$ appearing in the definition of $\Y$-symmetroidality can be weakly saturated.

\lemm \label{lemm--symmetroidal.Y.weakly.sat}
Let $S$, $\Y$, $\C$ be as in \refde{power.cofibration}.
If $S$ is $\Y$-symmetroidal, it is also $\ws \Y$-symmetroidal.
\xlemm
\pf
For a fixed $s \in S$, the functor $F_s\colon y \mapsto y \pp_{\Sigma_n} s^{\pp n}$ is cocontinuous.
In particular, $F_s(\ws \Y) \subset \ws {F_s(\Y)} \subset \ws \cof_\C = \cof_\C$ and likewise for acyclic $\Y$-symmetroidal maps.
\xpf

\defi
The cofibrant replacement of~$1$ in~$\Sigma_n^{\op,\proj}\V$ is denoted by~$\E\Sigma_n$.
(For $\V=\sSet$, this coincides with the usual definition of~$\E\Sigma_n$ as a weakly contractible simplicial set with a free $\Sigma_n$-action.)
\xdefi

\refpr{acyclic} is a key step in the proof of stability of symmetric h-monoidality and symmetroidality under left Bousfield localizations.
It relies on the following technical lemma.

\lemm\label{lemm--preparatory}
Suppose $\C$ is a symmetric monoidal, h-monoidal, flat model category,
$y\in\Sigma_n\C$ is any map, $s$ is a finite family of acyclic cofibrations with cofibrant domain that lies in some symmetric flat class~$S$,
and $y\pp s^{\pp n}$ is a weak equivalence in~$\C$ for some multiindex $n > 0$.
Then $y\pp_{\Sigma_n}s^{\pp n}$ is also a weak equivalence.
\xlemm

\pf
Let
$$\xymatrix{
A' \ar[r]^a_\sim \ar[d]_{y'} & A \ar[d]^y \\
B' \ar[r]_b^\sim & B
}$$
be the functorial cofibrant replacement of~$y \colon A \r B \in\Ar(\C)$ (in the projective model structure, so that $y'$~is a cofibration with a cofibrant domain).
Functoriality and the fact that $y\in\Ar(\Sigma_n\C)$ imply that $y'\in\Ar(\Sigma_n\C)$.
We claim that $y' \pp s^{\pp n}$ is a cofibrant replacement of $y \pp s^{\pp n}$ in~$\Ar(\C)$.
Let $t := s^{\pp n} \colon T \r S$.
The map $b \t S$ is a weak equivalence by the flatness assumption.
To see that $B' \t T \sqcup_{A' \t T} A' \t S \r B \t T \sqcup_{A \t T} A \t S$ is a weak equivalence we first note that these pushouts are homotopy pushouts by \refle{i.cofibrations}\refit{h.cofibration.containment} since $A \t t$ is an h-cofibration.
Thus, it suffices that the three individual terms in the pushouts are weakly equivalent, which again follows from flatness.
The claim is shown.
\ppar
Thus, we have
$$\hocolim_{\Sigma_n} (y \pp s^{\pp n}) = (\E\Sigma_n \t y' \pp s^{\pp n})_{\Sigma_n} \sim y \pp_{\Sigma_n} s^{\pp n}.$$
The last weak equivalence holds by symmetric flatness of $S$ since $\E\Sigma_n \t y' \r y' \r y$
is a weak equivalence by the unit axiom for the $\V$-enrichment
(note that the cofibrant replacement $\E\Sigma_n \r 1$ in $\Sigma_n^\proj \V$ is, in particular, a cofibrant replacement in~$\V$).
Finally, $y \pp s^{\pp n}$ is a weak equivalence in~$\C$ by assumption.
Therefore, the above homotopy colimit is a weak equivalence in $\C$.
\xpf

\prop \label{prop--acyclic}
The class of acyclic power cofibrations coincides with the intersection of $\we$ with the class of power cofibrations.

A $\Y$-symmetroidal class $S$ which consists of acyclic cofibrations with cofibrant source is acyclic $\Y$-\hskip0pt symmetroidal, provided that $\C$ is h-monoidal and flat and $S$ is symmetric flat in~$\C$.
\xprop

\pf
The first claim follows from the pushout product axiom.

For any $s \in S$ and any map $y \in \Y_n \subset \Mor (\Sigma_n \C)$, $y \pp s^{\pp n}$ is a weak equivalence in $\Sigma_n \C$ by assumption on the class $\Y$ (see \refde{symmetroidal}).
Now apply \refle{preparatory}.
\xpf

We now establish the compatibility of the three symmetricity properties with weak saturation.
Parts~\refit{symmetroidal.weakly.sat} and~\refit{power.weakly.sat} of \refth{symmetric.weakly.saturated}
are due to Gorchinskiy and Guletski\u\i~\cite[Theorem~5]{GorchinskiyGuletskii:Symmetric}.
Part~\refit{symmetric.flat.weakly.sat} extends arguments in \cite[Theorem~9]{GorchinskiyGuletskii:Positive},
which shows a weak saturation property for symmetrically cofibrant objects in a stable model category.
Of course, it also extends the analogous statement for nonsymmetric flatness (\refth{power.weakly.saturated}\refit{flat.weakly.sat}).
Likewise, \refit{symmetric.i.monoidal.weakly.sat}~extends the weak saturation property of h-cofibrations (see \refle{i.cofibrations}).
The proof of the closure under transfinite compositions in~\refit{symmetroidal.weakly.sat} is reminiscent of~\cite[\S4]{GorchinskiyGuletskii:Symmetric}.
See also the more recent accounts by White~\cite[Appendix~A]{White:Model} and Pereira~\cite[\S4.2]{Pereira:Cofibrancy}.
In the proof of the theorem, we will need a combinatorial lemma that we establish first.
Recall the conventions for multiindices in~\refde{multi}.

\lemm \label{lemm--combinatorial}
Let $X_0^{(i)} \To2{v_0^{(i)}}X_1^{(i)} \To2{v_1^{(i)}}X_2^{(i)}$, $1 \le i \le e$ be a finite family of composable maps in a symmetric monoidal category.
For a pair of multiindices $0\le k\le n$ of length $e$, we set $$m_k := \Sigma_n \cdot_{\Sigma_{n-k} \x \Sigma_k}  v_0^{\pp n-k} \pp v_1^{\pp k}.$$
\begin{enumerate}[(i)]
\item The map $(v_1 v_0)^{\pp n}\colon\bigppdom^n (v_1 v_0) \r X_2^{\t n}$ is the composition of pushouts (with the attaching maps constructed in the proof) of the maps~$m_k$ ($0\le k<n$),
and the map $m_n=v_1^{\pp n}$.
\item The map $\kappa\colon \bigppdom^n(v_1 v_0)\sqcup_{\bigppdom^n v_0}X_1^{\t n}\to X_2^{\t n}$
is the composition of pushouts of the maps~$m_k$ for $1\le k<n$, and the map~$m_n$.
(Here $1$ denotes the multiindex whose components are all equal to~1.)
\end{enumerate}
\xlemm

\pf
We interpret the composable pair~$(v_0,v_1)$ as a functor $v\colon3=\{0\to1\to2\}\to\C^I$, where $I = \{1, \ldots, e\}$.
Let $\mathcal E$ be the category of posets $C$ lying over $3^n=\prod_i3^{n_i}$ and let $\Sigma_n \mathcal E$ be those posets with a $\Sigma_n$-action which is compatible with the $\Sigma_n$-action on $3^n$.
For all posets considered below, the map to $3^n$ will be obvious from the context.
Consider the following functor:
$$Q\colon\Sigma_n \mathcal E \to\Sigma_n\C, \qquad(C \r 3^n)\mapsto\colim \left (C \To2{} 3^n \To2{v^n}\C^n \To2\otimes \C\right).$$
Being the composition of the two cocontinuous functors
$$\text{posets} / 3^n \To4{} \text{posets} / \C \To4\colim\C,$$
$Q$ is also cocontinuous.
The map $(v_1 v_0)^{\pp n}$ is obtained by applying $Q$ to the map
$$\iota\colon\{0,1,2\}^n \backslash \{1,2\}^n \r \{0,1,2\}^n,$$
which adds all tuples containing only $1$'s and $2$'s.
It is the composition of the maps
$$\iota_k\colon\{0,1,2\}^n\setminus\{1,2\}^n\cup\{\Sigma_n1^*2^{<k}\} \to \{0,1,2\}^n\setminus\{1,2\}^n\cup\{\Sigma_n1^*2^{\le k}\},$$
for $0\le k\le n$, with $\prod_i(n_i+1)$ maps in total.
The superscript~$*$ means that one adds as many elements as needed to get an $n$-multituple.
For multiindices the above statements should be interpreted separately for each component.
The map~$\iota_k$ adds the $\Sigma_n$-orbit~$O$ consisting of tuples with $k$~2's and $n-k$~1's,
i.e., $\Sigma_n1^{n-k}2^k$.
The cardinality of~$O$ is $n \choose k$.
For $o \in O$, consider the downward closure $D_{o}$ of $o$ and $C_o := D_o \backslash \{o\}$.

There is a pushout diagram in $\Sigma_n \mathcal E$
$$\vcenter{\xymatrix{
A := \coprod_{o \in O} C_o \ar[r]^(.35){\alpha_k} \ar[d]^{\mu_k} & \{0,1,2\}^n\setminus\{1,2\}^n\cup\{\Sigma_n1^*2^{<k}\} \ar[d]^{\iota_k}\cr
B := \coprod_{o \in O} D_o \ar[r] & \{0,1,2\}^n\setminus\{1,2\}^n\cup\{\Sigma_n1^*2^{\le k}\}.\cr}}\eqlabel{pushout.diagram}$$
(For $k=n$ the top horizontal row is an identity, so $\iota_n=\mu_n$ in this case.)
Any $o \in O$ determines a partition of $\coprod_in_i$ into $\coprod_i\{1\le j\le n_i\mid o_{i,j}=1\}$ and $\coprod_i\{1\le j\le n_i\mid o_{i,j}=2\}$.
Using this partition,
we have $D_o=\Sigma_{n-k}0^*1^*\times\Sigma_k0^*1^*2^*$ and
$C_o = \Sigma_{n-k}0^*1^{<n-k}\times\Sigma_k0^*1^*2^*\cup\Sigma_{n-k}0^*1^*\times\Sigma_k0^*1^*2^{<k}$.
Thus, the map $Q(C_o \r D_o)$ is just $v_0^{\pp n-k} \pp v_1^{\pp k}$.
Using the cocontinuity of $Q$, this shows $Q(\mu_k) = m_k$.

The second part now follows immediately from the above once we observe that the codomain of~$\iota_0$ is precisely the domain of the map under consideration.
\xpf

\theo \label{theo--symmetric.weakly.saturated}
Let $S$ be a class of morphisms in a symmetric monoidal model category $\C$.
We say some property of $S$ is \emph{stable under saturation} if it also holds for the weak saturation $\ws S$.
\begin{enumerate}[(i)]
\item \label{item--symmetric.flat.weakly.sat}
If $\C$ is pretty small and left proper, and $S$~is symmetric h-monoidal, then symmetric flatness of $S$ relative to a class $\Y = (\Y_n)$ of weak equivalences in $\Sigma_n \C$ is stable under saturation.
In particular, if some class of generating cofibrations in $\C$ is symmetric flat and symmetric h-monoidal, then $\C$ is symmetric flat.

\item
\label{item--symmetric.i.monoidal.weakly.sat}
If $\C$ is pretty small and left proper, then the property of being (acyclic) symmetric h-monoidal is stable under saturation.
In particular, if some class of generating (acyclic) cofibrations consists of (acyclic) symmetric h-cofibrations, then $\C$ is symmetric h-monoidal.

\item
\label{item--symmetroidal.weakly.sat}
Being $\Y$-symmetroidal (\refde{symmetroidal}) is stable under saturation.
In particular, if some class of generating (acyclic) cofibrations is (acyclic) $\Y$-sym\-met\-roi\-dal, then $\C$ is $\Y$-symmetroidal.

\item
\label{item--power.weakly.sat}
The same statement holds for power cofibrations.
\end{enumerate}
\xtheo

\pf
For a finite family of maps $v = (v^{(1)}, \ldots, v^{(e)})$ we use the multi-index notation of \refde{multi}.
We prove the statements by cellular induction, indicating the necessary arguments for each statement individually in each step.
The acyclic parts of~\refit{symmetric.i.monoidal.weakly.sat} and~\refit{symmetroidal.weakly.sat} are the same as the nonacyclic parts, so they will be omitted.
Fix an object $Y \in \Sigma_n \C$, respectively, a map $y \in \Y_n \subset \Mor \Sigma_n \C$.
For \refit{symmetric.flat.weakly.sat},~\refit{symmetroidal.weakly.sat} and~\refit{symmetric.i.monoidal.weakly.sat}, we write
$$g(v, n) := y \pp_{\Sigma_n} v^{\pp n}, \ \text{respectively, } g(v, n) := Y\otimes_{\Sigma_n} v^{\pp n}.$$
By \cite[Proposition~6.13]{Harper:Symmetric} or \refpr{pushout.product.pushout.morphisms}, $g(-,n)$ preserves pushout morphisms $\varphi\colon v \r v'$ (in the sense that, say, $\varphi^{(1)}$ is a pushout morphism and all other $\varphi^{(j)}$'s are identities) and retracts.
Thus, if $g(v, n)$ is an (acyclic) h-cofibration or (acyclic) cofibration, so is $g(v', n)$ by \refle{i.cofibrations}\refit{h.i.cofibrations.pushout} (which uses left properness).
This shows the stability of the properties of being symmetric h-monoidal and symmetroidal under cobase changes.
For~\refit{symmetric.flat.weakly.sat}, we additionally observe that $Y \t_{\Sigma_n} v^{\pp n}$ is an h-cofibration and similarly with $Y'$ since $S$ is symmetric h-monoidal by assumption.
By \refle{preparation.flat.weakly.sat} (more precisely, replace $\t$ there by $\t_{\Sigma_n}$), applied to $a=v^{\pp n}$ and $b=v'^{\pp n}$, we see that $g(v', n)$ is a weak equivalence since $g(v, n)$ is one.

We now show the stability of the three symmetricity properties relative to a class under transfinite composition: suppose $v^{(1)}$ is the transfinite composition
$$v^{(1)}\colon X^{(1)}_0 \To2{v^{(1)}_0}\cdots\to X^{(1)}_i \To2{v^{(1)}_i}X^{(1)}_{i+1}\to \cdots\to X^{(1)}_\infty = \colim X^{(1)}_i,$$
whose maps are obtained as pushouts
$$\xymatrix{A \ar[d] \ar[r]^{s \in S} \ar@{}[dr]|{(*)} & A' \ar[d]\cr
X := X^{(1)}_i \ar[r]^{x := v^{(1)}_i} & X' := X^{(1)}_{i+1}.\cr}$$
For the statements \refit{symmetric.flat.weakly.sat},~\refit{symmetric.i.monoidal.weakly.sat},~\refit{symmetroidal.weakly.sat}
we need to show that the map $g(v, n) = g((v^{(1)}, \ldots, v^{(e)}, n)$ is a weak equivalence, h-cofibration, or cofibration, respectively, provided that
$$\{v_i^{(1)}, i \le \infty, v^{(2)}, \ldots, v^{(e)}) \}$$
is a symmetric flat, symmetric h-monoidal, or symmetroidal class, respectively.
Applying this argument $e$ times gives the desired stability under transfinite compositions.
We write $r^{(1)}_i\colon X^{(1)}_0 \to X^{(1)}_i$ for the (finite) compositions of the $v^{(1)}_i$.
Consider
$$\id_{(X^{(1)}_0)^{\t n}} = (r^{(1)}_0)^{\pp n} \to (r^{(1)}_1)^{\pp n} \to \cdots \to (v^{(1)})^{\pp n}.$$
As an object of $\Sigma_n \Ar(\C)$,
$$g(v, n) = \colim_i g(\underbrace{(r_i^{(1)}, v^{(2)}, \ldots, v^{(e)})}_{=:v_i}, n) = \colim_i g(v_i, n),$$
since $-^{\pp n}$ preserves filtered colimits \cite[Corollary 4.4.5]{GambinoJoyal:Operads}.
We now show that $v_i$ is a symmetric flat (respectively, symmetric h-monoidal or symmetroidal) family,
so that $g(v_i)$ is a weak equivalence (respectively, h-cofibration or cofibration).
We consider the composition of two morphisms $r_0^{(1)}$ and $r_1^{(1)}$ only and leave the similar case of a finite composition of more than two maps to the reader.
By Lemmas \ref{lemm--pp.compositions} and~\ref{lemm--combinatorial}, $v_1^{\pp n}$ is the (finite) composition of pushouts of $\Sigma_n \cdot_{\Sigma_m} w^{\pp m}$, where
$w = (r_0^{(1)}, r_1^{(1)}, v^{(2)}, \ldots, v^{(e)})$,
and $m$ runs through multi-indices of length $e+1$ such that $0 \le m^{(1)} \le n^{(1)}$, $m^{(1)} + m^{(2)} = n^{(1)}$, and $m^{(k)} = n^{(k-1)}$ for $2 \le k \le e+1$.

For~\refit{symmetric.i.monoidal.weakly.sat}, each $g(w, m)  = y \pp_{\Sigma_m} w^{\pp m}$ is an h-cofibration.
Hence, so is $g(v_1, n)$ since h-cofibrations are stable under pushouts and (finite) compositions by \refle{i.cofibrations}.
By \refle{i.cofibrations}\refit{i.cofibrations.weakly.sat}, $g(v, n)$ is also an h-cofibration then.

Similarly, for~\refit{symmetroidal.weakly.sat}, each $g(w, m)$ is a cofibration, so that $g(v_1, n)$ is a cofibration.
By \refle{combinatorial}, $(v_1^{(1)} \circ v_0^{(1)})^{\pp n}$ is the composition of a pushout of $(v_0^{(1)})^{\pp n}$ and the map
$$\bigppdom^{n^{(1)}} (v_1^{(1)} \circ v_0^{(1)}) \sqcup_{\bigppdom^{n^{(1)}} (v_0^{(1)})} (X_1^{(1)})^{\t n} \r (X_2^{(1)})^{\t n}.$$
Here, as usual, $\bigppdom^{n^{(1)}} -$ denotes the domain of the $-^{\pp n^{(1)}}$.
The latter map is the composition of pushouts of the maps $g(w, m)$, where $w$ and $m$ are as above, except that now $0 \le m^{(1)} < n^{(1)}$.
Again, these are cofibrations, so the above map is a cofibration.
By \refle{sequential}\refit{cof.sequential}, $g(v, n)$ is therefore a cofibration.

For~\refit{symmetric.flat.weakly.sat}, each $g(w, m)$ is a weak equivalence.
The map $g(v_1, n)$ is the composition of pushouts of $g(w, m)$ along $Y \t_{\Sigma_n} \Sigma_n \cdot_{\Sigma_m} w^{\pp m} = Y \t_{\Sigma_m} w^{\pp m}$.
The latter map (and similarly for $Y'$) instead of $Y$ is an h-cofibration by the symmetric h-monoidality assumption.
Thus, the pushouts of $g(w, m)$, the compositions of which are $g(v_1, n)$, are weak equivalences by \refle{preparation.flat.weakly.sat} (again, replace $\t$ by $\t_{\Sigma_n}$ there).
We have shown that $g(v_1, n)$ is a weak equivalence.
By \refle{sequential}\refit{we.chain.colimits}, $g(v, n)$ is then also weak equivalence.

\refit{power.weakly.sat} can be shown using the same argument but considering $g(v) := v^{\pp n} \in \Sigma_n \C$ instead.
By \refre{symmetroidal} it is unnecessary to use multi-indices in this proof.
\xpf

\numberwithin{equation}{section}

\section{Transfer of model structures}
\label{chap--transfer}

In this section, we fix an adjunction
$$F:\C\leftrightarrows \D:G \eqlabel{adjunction}$$
such that $\C$ is a model category and $\D$ is complete and cocomplete.

\defi \label{defi--transfer}
A model structure on~$\D$ is {\it transferred along~$G$\/}
if the weak equivalences and fibrations in~$\D$ are those morphisms which are mapped by $G$ to weak equivalences and fibrations in $\C$, respectively.
\xdefi

If a transferred model structure on $\D$ exists, it is unique, so we also speak of \emph{the} transferred model structure.
In the sequel, we fix a Quillen adjunction as above such that the model structure on $\D$ is transferred from $\C$.

The next proposition describes basic properties of transferred model structures.
Part~\refit{transfer.left.proper}
is a special case of much more general left properness results by Batanin and Berger \cite[Theorem 2.11]{BataninBerger:Homotopy}.

\prop \label{prop--transfer.basic}
The following properties hold for a transferred model structure on~$\D$.
We write $I$ (respectively, $J$) for a class of generating (acyclic) cofibrations of $\C$.
\begin{enumerate}[(i)]
\item \label{item--transfer.cofibrations}
The class~$F(I)$ (respectively, $F(J)$) generates (acyclic) cofibrations of~$\D$.
In particular, if $\C$ is quasi-tractable, then so is~$\D$.
Moreover, if $\C$ is combinatorial  \cite[Definition~A.2.6.1]{Lurie:HTT}, then so is~$\D$, provided that $\D$ is locally presentable.
\item\label{item--transfer.pretty.small}
Suppose that $G$ preserves filtered colimits.
If $\C$ is pretty small, then so is~$\D$,
provided that $\D$ is locally presentable, or, more generally, $F(I')$ and $F(J')$ permit the small object argument,
where $I'$~and~$J'$ come from pretty smallness.
\item \label{item--transfer.left.proper}
Suppose that $G$ is cocontinuous and suppose that (a)~$G(F(I))$ consists of cofibrations or (b)~$\C$ is pretty small and $G(F(I))$ consists of h-cofibrations.
Then, if $\C$ is left proper, so is~$\D$.
\item \label{item--G.preserves.cofibrations}
If $G$ preserves filtered colimits and sends cobase changes of~$F(I)$ (respectively, cobase changes of~$F(I)$ along maps with cofibrant targets)
to cofibrations, then $G$ preserves cofibrations (respectively, cofibrations with cofibrant source).
\end{enumerate}
\xprop

\pf
\refit{transfer.cofibrations}:
Follows from $\inj \cof =\AF$ and $\inj \AC=\fib$.

\refit{transfer.pretty.small}:
By \refde{pretty.small}, there is another model structure~$\C'$ on the underlying category of~$\C$
with the same weak equivalences and a smaller class of cofibrations
that is generated by a set of morphisms with compact domains and codomains.
By assumption $F(\cof_{\C'})$ permits the small object argument and similarly for acyclic cofibrations.
This verifies the condition for the existence of the transfer of the model structure~$\C'$.
Thus, the model structure~$\C'$ transfers to a model structure~$\D'$ on the category underlying~$\D$ and its cofibrations are a subset of cofibrations of~$\D$.
The (co)domains of the generating set of cofibrations~$F(I')$ are compact because $G$ preserves filtered colimits and therefore $F$ preserves compact objects.

\refit{transfer.left.proper}:
By \refle{i.cofibrations}\refit{h.cofibration.containment} and \refle{i.cofibrations.detect} we have to show that $G(\cof_\D)$ consists of h-cofibrations in $\C$.
This follows from the assumptions (in case (b) use \refle{i.cofibrations}\refit{i.cofibrations.weakly.sat}).

\refit{G.preserves.cofibrations}:
Cofibrations in~$\D$ are retracts of transfinite compositions of cobase changes of elements in~$F(I)$.
All three operations are preserved by the functor~$G$ by assumption.
Thus, it is sufficient to observe that $G(F(I))$ consists of cofibrations in~$\C$, which are weakly saturated,
hence $G$ preserves cofibrations.
The preservation of cofibrations with cofibrant source is shown in the same way.
\xpf



From now on, we assume in addition that $\C$ and $\D$ are symmetric monoidal categories.
By \refpr{transfer.basic}\refit{transfer.cofibrations}, $\D$ is a monoidal model category provided that $F$ is strong monoidal.

\defi
\label{defi--Hopf.adjunction}
A {\it Hopf adjunction\/} is an adjunction between monoidal categories such that there is a functorial isomorphism
for $C\in\C$, $D\in\D$,
$$G(F(C)\otimes D)\cong C\otimes G(D).$$
\xdefi

\rema
If the monoidal products $\t_{\C}$ and $\t_\D$ are closed, this is equivalent to~$G$ being strong closed, i.e., internal homs are preserved up to a coherent isomorphism.
\xrema

\prop\label{prop--transfer.monoidal}
Suppose the adjunction $(F, G)$ is a Hopf adjunction.
Also suppose that $G$ preserves pushouts along maps of the form $D \t F(s)$, where $D \in \D$ is any object and $s$ is any morphism in $S$.
Let $S$ be a class of cofibrations in~$\C$.
We say that a property of the class~$S$ \emph{transfers}, if the same property holds for~$F(S)$.

\begin{enumerate}[(i)]
\item
\label{item--transfer.i.monoidal}
Suppose $\D$ is left proper.
Then the (acyclic) h-monoidality of $S$ transfers.
The h-monoidality of $\C$ transfers to~$\D$ if $\D$ is pretty small.

\item \label{item--transfer.flat}
The flatness of $S$ transfers.
The flatness of $\C$ transfers to $\D$ if $\D$ is pretty small and h-monoidal.

\item
\label{item--transfer.monoid.axiom}
If $G$ also preserves filtered colimits then the monoid axiom transfers from~$\C$ to~$\D$.
\end{enumerate}
\xprop

\pf
\refit{transfer.i.monoidal} and~\refit{transfer.flat} are shown exactly the same way as their symmetric counterparts, see Parts~\refit{transfer.symmetric.i.monoidal} and~\refit{transfer.symmetric.flat} of \refth{transfer.symmetric.monoidal}, using \refth{power.weakly.saturated} instead.


\refit{transfer.monoid.axiom}:
The preservation of colimits under $\t_{\D}$ and \refpr{transfer.basic}\refit{transfer.cofibrations},
the assumption that $G$ preserves the weak saturation, the Hopf adjunction property, and the monoid axiom for~$\C$ imply 
$$\eqalign{G(\ws{\D \t \AC_{\D}}) & \subset G(\ws{\D \t F(\AC_{\C})}) \subset \ws{G(\D \t F(\AC_{\C})}\cr
&= \ws{G(\D) \t \AC_{C}} \subset \ws{\C \t \AC_{\C}} \subset \we_\C.\cr}$$
\xpf

The following theorem shows that the three symmetricity properties interact well with transfers.
It is the symmetric counterpart of \refpr{transfer.monoidal}.

\theo
\label{theo--transfer.symmetric.monoidal}
Let
$F : \C \leftrightarrows \D : G$
be a Quillen adjunction of symmetric monoidal model categories such that the model structure on $\D$ is transferred from $\C$.
We assume $F$ is strong monoidal.
For Parts~\refit{transfer.symmetric.i.monoidal} and~\refit{transfer.symmetric.flat}, we also assume
\begin{enumerate}[(a)]
\item
\label{item--assumption.Hopf.adjunction}
the adjunction is a Hopf adjunction,

\item
\label{item--assumption.finite.colimits}
$G$ preserves finite colimits (including pushouts and $\Sigma_n$-coinvariants).
\end{enumerate}

Let $S$ be a class of cofibrations in $\C$.
We say that a property of the class $S$ \emph{transfers}, if the same property holds for $F(S)$.

\begin{enumerate}[(i)]
\item
\label{item--transfer.symmetric.i.monoidal}
Suppose $\D$ is left proper (a sufficient criterion is given in \refpr{transfer.basic}\refit{transfer.left.proper}).
Then the (acyclic) symmetric h-monoidality of $S$ transfers.
The symmetric h-monoidality of $\C$ transfers if, in addition, $\D$ is pretty small.

\item \label{item--transfer.symmetric.flat}
Symmetric flatness of $S$ transfers.
Moreover, the symmetric flatness of $\C$ transfers to $\D$ if, in addition, $\D$ is pretty small and symmetric h-monoidal.

\item \label{item--transfer.symmetroidal}
For some class $\Y$ of morphisms as in \refde{symmetroidal}, the $\Y$-symmetroidality of $S$ transfers in the sense that $\ws{F(S)}$ is $F(\Y)$-symmetroidal.
In particular, if $\C$ is $\Y$-symmetroidal, then $\D$ is $\ws{F(\Y)}$-symmetroidal.

\item \label{item--transfer.freely.powered}
Then the property of being freely powered transfers.
In particular, if $\C$ is freely powered, then so is $\D$.
\end{enumerate}
\xtheo

\pf
For all properties, the transfer for the given class $S$ is proven using a specific argument.
The transfer of the property from $\C$ to $\D$ follows from the fact that $F(\cof_\C)$ generates the cofibrations of $\D$ (\refpr{transfer.basic}\refit{transfer.cofibrations}), and likewise for acyclic cofibrations.
Then, a weak saturation property (indicated below) is used.
Let $s \in S$ be any map.

\refit{transfer.symmetric.i.monoidal}:
We need to show that $Y \t_{\Sigma_n} F(s)^{\pp n} = Y \t_{\Sigma_n} F(s^{\pp n})$ is an (acyclic) h-\hskip0pt cofibration for all $Y \in \Sigma_n\D$.
By \refle{i.cofibrations.detect} it is enough to show $G(Y \t_{\Sigma_n} F(s^{\pp n}))$ is an (acyclic) h-cofibration.
Using the Hopf adjunction property, strong monoidality of $F$ (which ensures that $F$ commutes with pushout products by \refpr{Ar.functorial}), and preservation of finite colimits under $G$, we compute this as $G(Y) \t_{\Sigma_n} s^{\pp n}$,
which indeed is an (acyclic) h-cofibration by the (acyclic) symmetric h-monoidality of~$S$.
The symmetric h-monoidality of $\C$ transfers to $\D$ by \refth{symmetric.weakly.saturated}\refit{symmetric.i.monoidal.weakly.sat}, using the left properness of $\D$.

\refit{transfer.symmetric.flat}:
For any weak equivalence $y$ in $\Sigma_n \D$ we have to show that $y \pp_{\Sigma_n} F(s)^{\pp n}$ is a weak equivalence.
As above, we have an isomorphism $G(y \pp_{\Sigma_n} F(s)^{\pp n})=G(y) \pp_{\Sigma_n}s^{\pp n}$.
This is a weak equivalence since $\C$~is symmetric flat.
To transfer the symmetric flatness of $\C$ we apply \refth{symmetric.weakly.saturated}\refit{symmetric.flat.weakly.sat} to $S = I_\C$, noting that $\D$ is h-monoidal by assumption and therefore left proper.

\refit{transfer.symmetroidal}:
As $F$ is strong monoidal and cocontinuous, $F(y) \pp_{\Sigma_n} F(s^{\pp n}) = F(y \pp_{\Sigma_n} s^{\pp n})$.
This shows the $F(\Y)$-symmetroidality since $F$ preserves cofibrations and acyclic cofibrations.
Then apply \refle{symmetroidal.Y.weakly.sat}.
The claim about the symmetroidality of $\D$ follows from \refth{symmetric.weakly.saturated}\refit{symmetroidal.weakly.sat}.

\refit{transfer.freely.powered}: Replace $y \pp_{\Sigma_n} s^{\pp n}$ by $s^{\pp n}$ in~\refit{transfer.symmetroidal} and use \refth{symmetric.weakly.saturated}\refit{power.weakly.sat}.
\xpf

\rema \label{rema--issue.symmetroidal}
If $\C$ is symmetroidal (i.e., symmetroidal with respect to the injective cofibrations in $\Sigma_n \C$), $\D$ need not be symmetroidal: for example, for $\C = \sSet$ and $\D = \Mod_R(\sSet)$ with $R=\Z/4$, i.e., simplicial sets with an action of $\Z/4$.
In this case, $R$ has a $\Z/2$-action, so $R$ is injectively cofibrant in $\Sigma_2 \Mod_R$, but $R \t_{R, \Sigma_2} R^{\t_R 2} = R / 2$ is not cofibrant as an $R$-module.
\xrema

We conclude this section by applying the criteria developed above to the case of the category
of modules over a commutative monoid~$R$ in a symmetric monoidal model category~$\C$.
An example of this situation occurs in the construction of unstable model structures on symmetric spectra,
which are by definition modules over a commutative monoid in symmetric sequences \cite[Theorem~5.1.2]{HoveyShipleySmith:Symmetric}.

As $R$ is commutative, the category~$\Mod_R$ of $R$-modules has a symmetric monoidal structure:
$$X\otimes_R Y:=\coeq(X\otimes R\otimes Y\rightrightarrows X\otimes Y).$$
The free-forgetful adjunction
$F = R \t - : \C \leftrightarrows \Mod_R : U$
has the following properties: $R \t -$ is strong monoidal since $(R\otimes X)\otimes_R(R\otimes Y)\cong R\otimes(X\otimes Y)$.
Moreover, it is a Hopf adjunction: $(R\otimes C)\otimes_R D\cong C\otimes D$.
Finally,~$U$ also has a right adjoint, the internal hom functor~$\Hom(R,-)$ (also known as the cofree $R$-module functor).
In particular, $U$ is cocontinuous.

The following theorem summarizes the properties of the transferred model structure on $\Mod_R$.
The existence of the model structure is due to Schwede and Shipley \cite[Theorem~4.1(2)]{SchwedeShipley:Algebras}.
As in \refth{transfer.symmetric.monoidal}, we say that some model-theoretic property \emph{transfers} if it holds for $\Mod_R$, provided that it does for $\C$.
The transfer of left properness to $\Mod_R$ (and much more general algebraic structures) was established by Batanin and Berger under the assumption that $\C$ is strong h-monoidal \cite[Theorems~2.11, 3.1b]{BataninBerger:Homotopy}.
The transfer of symmetric flatness, symmetric h-monoidality and symmetroidality is new.

\theo \label{theo--commutative.monoid}
Suppose $\C$ is a cofibrantly generated symmetric monoidal model category that satisfies the monoid axiom and $R$ is a commutative monoid in~$\C$.
The transferred model structure on~$\Mod_R$ exists and is a cofibrantly generated symmetric monoidal model category.

Combinatoriality, (quasi)tractability, admissible generation, pretty smallness, $\V$-enrichedness, and the property of being freely powered transfer from $\C$ to $\Mod_R$.
Moreover, if $\C$ is symmetroidal with respect to some class~$\Y$ (\refde{symmetroidal}),
then $\Mod_R$ is symmetroidal with respect to $\ws{R \t \Y}$, the weak saturation of maps of free $R$-module maps generated by all $y \in \Y$.

If either $R$ is a cofibrant object in~$\C$ or if $\C$ is pretty small and h-monoidal, then left properness transfers.

If $\C$ is pretty small and h-monoidal, then flatness, symmetric flatness, h-mon\-oidality, symmetric h-mon\-oidality, and the monoid axiom transfer from $\C$ to $\Mod_R$.
\xtheo

\pf
The existence of the transferred model structure follows from \cite[Theorem~11.3.2]{Hirschhorn:Model}
since $F(J)=R\otimes J$ and the class of $F(J)$-cellular maps consists of weak equivalences by the monoid axiom.
The transfer of combinatoriality, (quasi)\hskip0pt tractability, pretty smallness, enrichedness, and left properness were established in \refpr{transfer.basic}.
The transfer of flatness, h-monoidality, and the monoid axiom is shown in \refpr{transfer.monoidal}, while their symmetric counterparts are treated in \refth{transfer.symmetric.monoidal}.
(We don't need an additional left properness assumption on $\C$ since any h-monoidal model category is left proper \cite[Lemma 1.8]{BataninBerger:Homotopy}.)
\xpf

\section{Left Bousfield localization}
\label{chap--localization}

Left Bousfield localizations of various types (e.g., ordinary, enriched, monoidal) of model categories
present reflective localizations of the corresponding locally presentable $\infty$-categories,
i.e., they invert the reflective saturation of a given class of maps in a (homotopy) universal fashion.
If the Bousfield localization of a given model category exists,
it can be constructed as a model structure on the same underlying category, with a larger class of weak equivalences
and the same class of cofibrations.
Examples for left Bousfield localizations abound, e.g.,
local model structures on simplicial presheaves (see \refchap{examples}) and the stable model structure on symmetric spectra are left Bousfield localizations.
In this section, we carry h-monoidality and flatness and their symmetric counterparts along a Bousfield localization.
An example application in the context of symmetric spectra is given in \csp{\refsect{stable}}.
The idea of combining h-monoidality and flatness was independently used by White~\cite{White:Monoidal}.

\defi \label{defi--localization}
\cite[Theorem~3.3.19]{Hirschhorn:Model}, \cite[Definitions~4.2, 4.42]{Barwick:Left}
Let $W$ be either of the following bicategories (specified by their objects, 1-morphisms, and 2-morphisms):
\begin{enumerate}[(a)]
\item model categories, left Quillen functors, and natural transformations;
\item
\label{item--monoidal.localization}
(symmetric) monoidal model categories, strong (symmetric) monoidal left Quillen functors, and (symmetric) monoidal natural transformations;
\end{enumerate}
Suppose $\C\in W$ and $S$ is a class of morphisms in~$\C$.
A {\it left Bousfield localization of~$\C$ with respect to~$S$\/}
is a 1-morphism $j\colon\C\to\BL_S\C$
such that precomposition with~$j$ induces an equality between
the category of morphisms $\BL_S\C\to\mathcal E$ (which are, in particular, left Quillen functors)
and the category of morphisms $\C\to\mathcal E$ whose left derived functors
send elements of~$S$ to weak equivalences in~$\mathcal E$.
\xdefi

For Part~\refit{monoidal.localization}, we write $\BL^\t$ instead of~$\BL$ and refer to this as a \emph{monoidal Bousfield localization} (the terminology is due to White~\cite{White:Monoidal}).
By \cite[Lemma~26]{GorchinskiyGuletskii:Symmetric}, the underlying model category of a monoidal localization is given by $U(\BL^\t_S\C)=\BL_{S^\otimes}U(\C)$,
where $S^\otimes$ is the monoidal saturation of~$S$,
which consists of the derived monoidal products of the elements of~$S$ and the objects of~$\C$
(or some class of homotopy generators of~$\C$, e.g., the set of domains and codomains of some set of generating cofibrations of~$\C$).
Fibrant objects in $\BL_{S^\t} \C$ are those fibrant objects $W$ in $\C$ such that the derived internal hom
(as opposed to the derived mapping space, which appears in nonmonoidal localizations)
$\hHom_\C(\xi, W)$ is a weak equivalence in $\C$ for any $\xi \in S$ \cite[4.46.4]{Barwick:Left}.

\rema
The above definition talks about equality of categories to ensure that the underlying category of a left Bousfield localization does not change.
One can replace equality with isomorphism or equivalence, which would yield an isomorphic or equivalent underlying category.
\xrema

\rema \label{rema--monoidal.enriched}
The above definition admits an obvious variant for $\V$-enriched localizations.
If $\C$ is $\V$-enriched and monoidal,
then monoidal localizations and $\V$-enriched monoidal localizations agree,
which immediately follows from the description of the monoidal saturations above and its enriched analog.
\xrema

In the following two theorems,
we say that a property of a class $S$ of cofibrations in $\C$ \emph{localizes} if it holds for $S$ regarded as a class of cofibrations in~$\D := \BL^\otimes_T\C$.
Likewise, we say that some property of~$\C$ localizes, if it also holds for~$\D$.

\prop \label{prop--monoidal.localization}
Let $\C$ be a (symmetric) monoidal model category such that the monoidal left Bousfield localization $\D :=\BL^\otimes_T\C$ with respect to some class $T$ exists.

\begin{enumerate}[(i)]
\item \label{item--localization.left.proper} \cite[Proposition 3.4.4]{Hirschhorn:Model}
Left properness of $\C$ localizes.
\item \label{item--localization.flat}
Flatness of~$S$ localizes.
In particular, the flatness of~$\C$ localizes.
\item \label{item--localization.i.cofibrations}
If $\C$ is left proper, any (acyclic) h-cofibration $f$ in~$\C$ is also an (acyclic) h-cofibration in~$\D$.
\item \label{item--localization.i.monoidal}
If $\C$ is left proper, quasi-tractable, pretty small, and flat, then the h-monoidality of $S$ or of $\C$ localizes.
\item \label{item--localization.monoid.axiom}
Pretty smallness localizes.
If $\D$ is pretty small and h-monoidal, then $\D$ satisfies the monoid axiom.
\end{enumerate}
\xprop

\pf
\refit{localization.flat}: We have to show that $y \pp s$ is a weak equivalence in $\D$ for all weak equivalences $y$ in $\D$ and $s \in S$.
By the pushout product axiom (of $\D$), we may assume $y$ is a trivial fibration in $\D$ or, equivalently, one in $\C$.
Now invoke the flatness of $S$ in $\C$ and use $\we_\C \subset \we_\D$.

\refit{localization.i.cofibrations}:
The acyclic part follows from the nonacyclic one and the inclusion $\we_\C \subset \we_\D$.
Given a diagram $A \leftarrow B \To1f C$, where $f$ is an h-cofibration in $\C$, we have to show by~\refit{localization.left.proper} and \refle{i.cofibrations}\refit{h.cofibration.containment} that $C \sqcup_B A$ is a homotopy pushout in~$\D$.
The identity functor $\Fun (\bullet \leftarrow \bullet \r \bullet, \C) \r \Fun (\bullet \leftarrow \bullet \r \bullet, \D)$ is a left Quillen functor if we equip both functor categories with the projective model structure.
Since it also preserves all weak equivalences, it preserves homotopy colimits, i.e., sends the homotopy pushout $C\sqcup_B A\sim C\sqcup_B^{h,\C}A$ to a homotopy pushout in~$\D$.

\refit{localization.i.monoidal}: As the cofibrations in $\C$ and $\D$ are the same, the nonacyclic part of the h-monoidality of $\D$ follows from~\refit{localization.i.cofibrations}.
The acyclic part of h-monoidality of $\D$ now follows from \refit{localization.left.proper} and \refco{acyclic.h.monoidal}.

\refit{localization.monoid.axiom}:
The first claim is clear since $\cof_\C = \cof_\D$.
The second is \refle{i.monoidal.monoid.axiom} again.
\xpf

\theo \label{theo--symmetric.monoidal.localization}
Let $\C$ be a (symmetric) monoidal model category such that the monoidal left Bousfield localization $\D :=\BL^\otimes_T\C$ with respect to some class $T$ exists.
\begin{enumerate}[(i)]
\item \label{item--localization.symmetric.flat}
Let $\Y = (\Y_n)$ be some classes of morphisms in $\Sigma_n \C$.
The property of $S$ of being symmetric flat with respect to $\Y$ localizes.
Moreover, the symmetric flatness of $S$ and of $\C$ localizes.
\item \label{item--localization.symmetric.i.monoidal}
If $\C$ is left proper and $\D$ is quasi-tractable, pretty small and symmetric flat, then the symmetric h-monoidality of $S$ or of $\C$ localizes.
\item \label{item--localization.symmetroidal}
The property of $S$ of being (acyclic) $\Y$-symmetroidal localizes provided that $\D$ is flat and h-monoidal and provided that $S$ consists of cofibrations with cofibrant source and is symmetric flat in~$\D$.
In particular, if $\D$ is h-monoidal and symmetric flat and $\C$ is $\Y$-symmetroidal then $\D$ is also $\Y$-symmetroidal.
\item \label{item--localization.freely.powered}
The property of being freely powered localizes.
\end{enumerate}
\xtheo

\pf
\refit{localization.symmetric.flat}:
The first claim is obvious from $\we_\C \subset \we_\D$.
For the second claim we have to show that $z \pp_{\Sigma_n} s^{\pp n}$ is a weak equivalence in $\C$ for all $z \in \we_{\Sigma_n^\proj \D}$ and $s \in S$.
We write $z = y \circ c$, where $y \in \AF_{\Sigma_n^\proj \D} = \AF_{\Sigma_n^\proj \C}$, and $c \in \AC_{\Sigma_n^\proj \D}$.
The map $z \pp_{\Sigma_n} s^{\pp n}$ is the composition of a pushout of $c \pp_{\Sigma_n} s^{\pp n}$, followed by $y \pp_{\Sigma_n} s^{\pp n}$.
The latter is a weak equivalence since $S$ is symmetric flat with respect to $\we_{\Sigma_n^\proj \C}$ by assumption.
The former is an acyclic cofibration: for this we may by cocontinuity assume $c = \Sigma_n \cdot h$ is a generating trivial cofibration, which yields
$h \pp s^{\pp n}$, itself an acyclic cofibration in~$\C$ by the pushout product axiom.

\refit{localization.symmetric.i.monoidal}:
As (acyclic) h-cofibrations of $\C$ are contained in the ones of $\D$ (\refpr{monoidal.localization}\refit{localization.i.cofibrations}), a class $S$ which is (acyclic) symmetric h-monoidal in $\C$ is also (acyclic) symmetric h-monoidal in $\D$.

Now suppose that $\C$ is symmetric h-monoidal.
We want to show that (acyclic) $\D$-cofibrations form an (acyclic) symmetric h-monoidal class (in $\D$).
Again using the above fact, it is enough to show the acyclic part.
Once again, we may restrict to generating acyclic cofibrations (\ref{theo--symmetric.weakly.saturated}\refit{symmetric.i.monoidal.weakly.sat}).
Thus, let $s$ be a finite family of generating acyclic cofibrations in $\D$.
By quasi-tractability, we may assume they have cofibrant domains.
Setting $y \colon \emptyset \r Y$, the pushout product $y \pp s^{\pp n}$ is just $Y \t s^{\pp n}$,
which is a weak equivalence by the h-monoidality of $\D$ ensured by \refpr{monoidal.localization}\refit{localization.i.monoidal}.
Using the flatness and h-monoidality of $\D$ (\refpr{monoidal.localization}\refit{localization.flat},~\refit{localization.i.monoidal}), \refle{preparatory} applies to $s$ and $y$ and shows that $Y \t_{\Sigma_n} s^{\pp n}$ is a weak equivalence.

\refit{localization.symmetroidal}: The stability of the nonacyclic part of $\Y$-symmetroidality is obvious.
The acyclic part follows from \refpr{acyclic}, using the cofibrancy assumption and the symmetric flatness of $S$ in $\D$.
Similarly, by \ref{theo--symmetric.weakly.saturated}\refit{symmetroidal.weakly.sat}, the symmetroidality of $\D$ follows by using a set $S$ of generating acyclic cofibrations (of $\D$) with cofibrant domain, which is possible thanks to the tractability of~$\D$.

\refit{localization.freely.powered}: This follows from \refpr{acyclic}.
\xpf

\numberwithin{equation}{subsection}

\section{Examples of model categories}
\label{chap--examples}
\label{chap--examples.symmetric}

We discuss the model-theoretic properties of \refchap{model}, \refsect{i.monoidal.flat}, and \refchap{symmetric} for simplicial sets, simplicial presheaves, simplicial modules, topological spaces, chain complexes, and symmetric spectra.

\subsection{Simplicial sets}
\label{sect--simplicial.sets}

The most basic example of a monoidal model category is the category $\sSet$ of simplicial sets equipped
with the cartesian monoidal structure $A \t B = A \x B$ and the Quillen model structure, see, e.g., \cite[Theorem~I.11.3]{GoerssJardine:Simplicial}.
All objects are cofibrant, so $\sSet$ is left proper, flat, and h-monoidal.

Simplicial sets are symmetroidal:
given any monomorphism $y \in \Sigma_n \sSet$ and a finite family of monomorphisms $v \in \sSet$, $y \pp_{\Sigma_n} v^{\pp n}$ is a monomorphism.
Indeed, $y \pp v^{\pp n}$ is a $\Sigma_n$-equivariant monomorphism and passing to $\Sigma_n$-orbits preserves monomorphisms.
By \refth{symmetric.weakly.saturated}\refit{symmetroidal.weakly.sat}, the acyclic part of symmetroidality follows
if $y \pp_{\Sigma_n} v^{\pp n}$ is a weak equivalence for any~$y$ as above and any finite family of horn inclusions $v \colon \Lambda^m_k \r \Delta^m$ (where $m$ and $k$ are multiindices).
To this end we first construct a homotopy $h \colon \Lambda\times\Delta^m\to\Delta^m$ from the identity map $\Delta^m\to\Delta^m$ to the composition $\Delta^m\To1{}\Delta^0\To1k\Delta^m$
such that $\Lambda^m_k\subset\Delta^m$ is preserved by the homotopy.
Here $\Lambda$ is the 2-horn, which can be depicted as $0 \r 1 \gets 2$.
We parametrize $h$ by~$\Lambda$ and not by the usual~$\Delta^1$ since~$\Delta^m$ is not fibrant.
The map $h$ is uniquely specified by its value on vertices, i.e., $\{0,1,2\}\times\{0,\ldots,m\}\to\{0,\ldots,m\}$.
We have $(0,i)\mapsto i$, $(1,i)\mapsto\max(k,i)$, $(2,i)\mapsto k$.
Thus, we have constructed a simplicial deformation retraction $\Lambda\times(\Lambda^m_k\to\Delta^m)\to(\Lambda^m_k\to\Delta^m)$
that contracts the inclusion $\Lambda^m_k\to\Delta^m$ to the identity map $\Delta^0\to\Delta^0$.
(Morphisms of maps are commutative squares, as usual.)
The map $h$ gives rise to a simplicial deformation retraction
$$\Lambda\times(y\pp_{\Sigma_n}v^{\pp n}) \stackrel \Delta \to (\Lambda^{\x n} \times(y\pp v^{\pp n}))_{\Sigma_n} \cong y\pp_{\Sigma_n}(\Lambda\times v)^{\pp n} \stackrel h \to y\pp_{\Sigma_n}v^{\pp n}$$
using the fact that the diagonal $\Delta\colon \Lambda \r \Lambda^{\x n}$ is $\Sigma_n$-equivariant.
It contracts the map $y\pp_{\Sigma_n}v^{\pp n}$ to the map $y\pp_{\Sigma_n}(\id_{\Delta_0})^{\pp n}$.
For $n>0$ the latter map is the identity map on the codomain of~$y$, in particular, a weak equivalence,
hence so is $y\pp_{\Sigma_n}v^{\pp n}$.

Symmetroidality and cofibrancy of all objects implies that $\sSet$ is symmetric h-monoidal.
The category $\sSet$ is far from freely powered: the map $(\partial \Delta^1 \r \Delta^1)^{\pp 2}$ is not a $\Sigma_2$-projective cofibration, since $\Sigma_2$ does not act freely on the complement of the image.
Simplicial sets are not symmetric flat:  $\E\Sigma_n \r *$ is $\Sigma_n$-equivariant and a weak equivalence of the underlying simplicial sets,
but $\mathrm B \Sigma_n := (\E\Sigma_n)_{\Sigma_n} \r *$ is not a weak equivalence: recall that $\mathrm B \Sigma_2$ is weakly equivalent to~$\RR\P^\infty$,
the infinite real projective space.

Similar statements hold for pointed simplicial sets equipped with the smash product.

The category $\sSet$ also carries the Joyal model structure \cite[Theorem 2.2.5.1]{Lurie:HTT}.
It is an interesting question whether it is symmetric h-monoidal.

\subsection{Simplicial presheaves}
\label{sect--simplicial.presheaves}

A more general example  is the category $\sPSh(S)=\Fun(S^\op,\sSet)$ of simplicial presheaves on some site~$S$.
The \emph{projective model structure} on this category is transferred from the Quillen model structure on $\sSet$ along
$$\prod_{X \in S} \sSet \leftrightarrows \sPSh(S).\eqlabel{transfer.presheaf}$$
It is pretty small by \ref{prop--transfer.basic}\refit{transfer.pretty.small} and left proper by \ref{prop--transfer.basic}\refit{transfer.left.proper}.
The monoid axiom, h-\hskip0pt monoidality, flatness, and symmetric h-monoidality
follow from the corresponding properties of the injective model structure (note that the class of h-cofibrations only depends on the weak equivalences).
Alternatively, even though \refeq{transfer.presheaf}~is not a Hopf adjunction, the arguments of \refpr{transfer.monoidal} can be generalized to~\refeq{transfer.presheaf}.
The projective model structure is not in general symmetroidal (for $X \in S$, $(X^n)_{\Sigma_n}$ is in general not projectively cofibrant).

In the \emph{injective model structure} on $\sPSh(S)$, weak equivalences and cofibrations are checked pointwise.
It is combinatorial \cite[Proposition~A.2.8.2]{Lurie:HTT} and therefore tractable.
It is pretty small (as the second model structure in \refde{pretty.small}, take the projective structure).
Since all objects are cofibrant, it is left proper, h-monoidal and flat.
The symmetric monoidality, symmetric h-monoidality and symmetroidality (with respect to injective cofibrations $\Y_n = \cof_{\Sigma_n^\inje \sPSh(S)}$) follows from the one of $\sSet$.

There are various \emph{intermediate model structures} on $\sPSh(S)$ (see \cite{Jardine:Intermediate}),
such as Isaksen's flasque model structure~\cite{Isaksen:Flasque}.
They also have pointwise weak equivalences but other choices of cofibrations which lie between projective and injective cofibrations.
For such intermediate model structures, monoidality, h-monoidality, symmetric h-monoidality, symmetroidality, the monoid axiom, and flatness follow from the injective case and pretty smallness follows from the projective case.

The properties mentioned above are stable under \emph{Bousfield localization}.
For example, given some Grothen\-dieck topology $\tau$ on the site $S$, the $\tau$-local projective model structure is the left Bousfield localization of the projective model structure with respect to $\tau$-hypercovers \cite[Theorem~6.2]{DuggerHollanderIsaksen:Hypercovers}.
Since hypercovers are stable under product with any $X \in S$ by \cite[Proposition~3.1]{DuggerHollanderIsaksen:Hypercovers}, this is a monoidal localization.
It is also $\sSet$-enriched by \cite[Theorem~4.1.1(4)]{Hirschhorn:Model}.
By \refpr{monoidal.localization}, the localized model structure is again left proper, tractable, monoidal and h-monoidal, pretty small, flat, and satisfies the monoid axiom.
It is symmetric h-monoidal at least if $\tau$ has enough points, for in this case local weak equivalences are maps which are stalkwise weak equivalences \cite[page~39]{Jardine:Simplicial}.

\subsection{Simplicial modules}
\label{sect--simplicial.modules}

Let $R$ be a commutative simplicial ring and consider the transferred model structure on simplicial $R$-modules via the free-forgetful adjunction
$$R[-] : \sSet \rightleftarrows \sMod_R : U.$$
The model category $\sMod_R$ is pretty small by \refpr{transfer.basic}.
As for chain complexes, $\sMod_R$ is flat, but not symmetric flat (unless $R$ is a rational algebra).

Simplicial $R$-modules are symmetric h-monoidal.
The nonacyclic part follows from the fact that monomorphisms, i.e., \emph{injective cofibrations}, of simplicial $R$-modules are h-cofibrations.

We reduce the acyclic part of symmetric h-monoidality of $\sMod_R$ to the one of $\sSet$
via the cocontinuous strong monoidal functor $R[-]\colon(\sSet,\times)\to(\sMod_R,\otimes)$, which preserves weak equivalences.
Given any object $Y\in\Sigma_n\sMod_R$
and any finite family~$w$ of generating cofibrations of~$\sMod_R$, i.e., $w=R[v]$,
we have a deformation retraction
$$\eqalign{R[\Lambda]\otimes(Y \t_{\Sigma_n} R[v]^{\pp n}) & \To3{R[\Delta]} (R[\Lambda]^{\t n} \otimes Y \t_{\Sigma_n}R[v]^{\pp n})_{\Sigma_n}\cr
&\cong Y \t_{\Sigma_n}(R[\Lambda \x v])^{\pp n} \To3{R[h]} Y \t_{\Sigma_n}R[v]^{\pp n}\cr}$$
of $Y\t_{\Sigma_n}w^{\pp n}$ to a weak equivalence,
which shows that the former is also a weak equivalence.

Simplicial $R$-modules are symmetroidal with respect to the class $\{R[\cof_{\Sigma_n^\inje \sSet}]\}$,
which follows immediately from the symmetroidality of simplicial sets and cocontinuity and strong monoidality of~$R[-]$.
Note that $\sMod_R$ is not symmetroidal, as can be shown as in \refre{issue.symmetroidal}.

\subsection{Chain complexes}
\label{sect--chain.complexes}

The category $\Ch_R := \Ch (\Mod_R)$ of unbounded chain complexes of $R$-modules, for some commutative ring~$R$, carries the \emph{projective model structure}
whose weak equivalences are the quasiisomorphisms and fibrations are the degreewise epimorphisms.
It is enriched over $\Ch (\Mod_\Z)$  (equipped with the projective model structure).
The generating (acyclic) cofibrations are given by all shifts of the canonical inclusion $[0 \r R] \r [R \To2\id R]$ ($[0 \r 0] \r [R \To2\id R]$, respectively)
\cite[Definition~2.3.3, Theorem~2.3.11]{Hovey:Model}.
In particular, the model structure is tractable and pretty small.
It is flat, as can be seen using \refth{power.weakly.saturated}\refit{flat.weakly.sat}.
The category is h-monoidal by \cite[Corollary 1.14]{BataninBerger:Homotopy}.

It is \emph{not} symmetric flat, for the same reason as $\sSet$ above.
Moreover, it is neither symmetric h-monoidal nor symmetroidal: for the chain complex $A = [\Z \To2\id \Z]$ in degrees 1 and 0, we have
$$A^{\t 2} = [\Z \To5{(1, -1)} \Z \oplus \Z \To5+ \Z],$$
where from left to right we have the sign representation, the regular and the trivial representation of  $\Sigma_2$.
However, $(A^{\t 2})_{\Sigma_2} = [\Z/2 \To2{} \Z \To2\id \Z]$ is not exact nor cofibrant.

By the Dold-Kan correspondence $N\colon (\sMod_R, \x) \rightleftarrows (\Ch^+_R, \t)$ between simplicial $R$-modules and connective chain complexes of $R$-modules, the projective model structures correspond to each other.
However, $N$ fails to be a strong symmetric monoidal functor.
Instead, $\x$ corresponds to the {\it shuffle\/} tensor product~$\widetilde \t$ of chain complexes, which is much bigger than the usual tensor product.
According to~\refsect{simplicial.modules}, $(\Ch^+_R, \widetilde \t)$ is symmetric h-monoidal.
The reason why a similar argument fails for~$\otimes$ is that the (smaller) ordinary tensor product fails to allow for a $\Sigma_n$-equivariant diagonal map for an interval object.

If $R$ contains $\Q$, the picture changes drastically: every $R$-module $M$ with a $\Sigma_n$-action is projective as an $R$-module if and only if it is projective as an $R[\Sigma_n]$-module (Maschke's theorem).
Thus, the projective and injective model structure (with respect to the $\Sigma_n$-action) on $\Sigma_n \Ch(\Mod_R)$ agree.
Therefore, $\Ch (\Mod_R)$ is symmetric flat and freely powered (and therefore symmetroidal and symmetric h-monoidal).

With appropriate additional assumptions, the statements above can be generalized to chain complexes in a Grothendieck abelian category $\cA$.
For example, flatness and h-monoidality of $\Ch (\cA)$ require that projective objects $P \in \cA$ are flat, i.e., $P\t-$ is an exact functor.

\subsection{Topological spaces}
\label{sect--topological.spaces}

The category $\Top$ of compactly generated weakly Hausdorff topological spaces carries the Quillen model structure
which is transferred from $\sSet$ via the singular simplicial set functor.
This model category is left proper \cite[Theorem 13.1.10]{Hirschhorn:Model}, monoidal \cite[Corollary~4.2.12]{Hovey:Model}, and h-monoidal \cite[Example~1.15]{BataninBerger:Homotopy}.
It is cellular \cite[Propositions~4.1.4]{Hirschhorn:Model}, though not locally presentable and therefore not combinatorial.
However, it is admissibly generated.
Alternatively, one can use Smith's $\Delta$-generated topological spaces, which are combinatorial.

Topological spaces are \emph{not} pretty small.
However, since closed inclusions are stable under $\Sigma_n$-coinvariants, products with arbitrary spaces and pushout products, and compact spaces are compact relative to closed inclusions \cite[Lemma 2.4.1]{Hovey:Model}, they are strongly admissibly generated in the sense below.

Recall from the definition of smallness of an object $A \in \C$ relative to a subcategory $\D \subset \C$ from \cite[Definition 10.4.1]{Hirschhorn:Model} or \cite[Definition 2.1.3]{Hovey:Model}.
By definition, any object in a combinatorial model category is small with respect to $\C$, so $\C$ is automatically admissibly generated in the sense below.
Topological spaces are not combinatorial, but strongly admissibly generated by the above remarks and \refpr{admissibly.generated.robust}\refit{admissibly.generated.weakly.sat}.

\defi
\label{defi--admissibly.generated}
\label{defi--strongly.admissibly.generated}
A cofibrantly generated symmetric monoidal model category~$\C$ is \emph{admissibly generated} (respectively, \emph{strongly admissibly generated}) relative to a class~$S$ of morphisms in~$\C$ if all (co)domains of a set of generating cofibrations are small (respectively, compact) with respect to the subcategory
$$\cell (Y\otimes_{\Sigma_n} s^{\pp n}) \eqlabel{admissibly.generated}$$
for any finite family~$s \subset S$, any multi-index $n > 0$, and any object~$Y \in \Sigma_n \C$.
We call~$\C$ \emph{(strongly) admissibly generated} if it is (strongly) admissibly generated relative to the cofibrations~$\cof_\C$.
\xdefi


For the purposes of this paper and also \cite{PavlovScholbach:Operads}, \cite{PavlovScholbach:Spectra}, strongly admissibly generated model categories are just as good as pretty small combinatorial model categories, as shown by the following results.
\refpr{admissibly.generated.robust} shows that (strong) admissible generation is easy to establish in practice and robust under standard operations on model categories.

\prop
\label{prop--strongly.admissibly.generated}
Suppose $\C$ is a symmetric monoidal model category.
\begin{enumerate}[(i)]
\item
\label{item--1}
(Weak saturation of (symmetric) h-monoidality and (symmetric) flatness.)
\refth{power.weakly.saturated}\refit{flat.weakly.sat},~\refit{i.monoidal.weakly.sat},
\refth{symmetric.weakly.saturated}\refit{symmetric.flat.weakly.sat},~\refit{symmetric.i.monoidal.weakly.sat} continue to hold if we replace ``pretty small'' in these statements by ``strongly admissibly generated''.
\item
\label{item--2}
(Transfer and localization of symmetric flatness and symmetric h-monoidality.)
\refpr{transfer.monoidal}, 
\refth{transfer.symmetric.monoidal}, 
\refpr{monoidal.localization}, and
\refth{symmetric.monoidal.localization}, 
continue to hold if we replace ``pretty small'' by ``strongly admissibly generated''.
\item
\label{item--3}
(Transfer of left properness.)
\refpr{transfer.basic}\refit{transfer.left.proper} continues to hold if we replace ``pretty small'' by the condition that (co)domains of a set $I$ of generating cofibrations of $\C$ are compact relative to pushouts of maps $F(I)$.
\end{enumerate}
\xprop

\pf
Analogously to \refle{sequential}\refit{we.chain.colimits}, a filtered colimit $f_\infty$ of weak equivalences $f_i$ is a weak equivalence,
provided that (co)domains of the generating cofibrations of $\C$ are compact relative to the class spanned by the acyclic cofibrations and the transition maps $x_i$, $y_i$.
Similarly, if this size condition is satisfied, $f_\infty$ is an h-cofibration provided that the $f_i$ and the maps $X_{i+1} \sqcup_{X_i} Y_i \r Y_{i+1}$ are h-cofibrations.
This refines \refle{i.cofibrations}\refit{i.cofibrations.weakly.sat}.

\refit{1}: To show the weak saturation of symmetric h-monoidality as in \ref{theo--symmetric.weakly.saturated}\refit{symmetric.i.monoidal.weakly.sat} using only that $\C$ is strongly admissibly generated, we use (cf.~the proof of~\ref{theo--symmetric.weakly.saturated}\refit{symmetroidal.weakly.sat})
that the transition maps appearing in the proof of \ref{theo--symmetric.weakly.saturated}\refit{symmetric.i.monoidal.weakly.sat}
are precisely of the form as in~\refeq{admissibly.generated}.

As for the stability of symmetric flatness under weak saturation (\ref{theo--symmetric.weakly.saturated}\refit{symmetric.flat.weakly.sat}),
it is enough to show that for a transfinite composition~$s$ of symmetric flat maps~$s_j$, and a weak equivalence~$y$,
the filtered colimit $y \pp_{\Sigma_n} s^{\pp n} = \colim_i y \pp_{\Sigma_n} t_i^{\pp n}$ is also a filtered homotopy colimit,
where $t_i = s_i \circ \cdots \circ s_0$ are the (finite) compositions of $s_j$.
By the above variant of \refle{sequential}\refit{we.chain.colimits}, this is true if the (co)domains of generating cofibrations are compact relative to the transition maps of this filtered colimit.
By \refle{combinatorial} and its proof, especially~\refeq{pushout.diagram}, these transition maps are given by $y \pp_{\Sigma_n} Q(\alpha_k)$, so this is true again since $\C$ is strongly admissibly generated.

\refit{2}: The indicated statements use pretty smallness only to invoke \refth{symmetric.weakly.saturated}.

\refit{3}: This follows from the above variant of \refle{i.cofibrations} and the proof of \refpr{transfer.basic}\refit{transfer.left.proper}.
\xpf

\prop
\label{prop--admissibly.generated.robust}
Suppose $\C$ is a symmetric monoidal model category.
\begin{enumerate}[(i)]
\item
\label{item--admissibly.generated.weakly.sat}
As in \refth{symmetric.weakly.saturated}, the property of being admissibly generated relative to~$S$ is stable under saturation.
Therefore, if $\C$ is cofibrantly generated and admissibly generated relative to some set of generating cofibrations, it is admissibly generated.

\item
\label{item--transfer.admissibly.generated}
In the situation of \refth{transfer.symmetric.monoidal} (including the assumptions (a)--(b) there), suppose $G$ preserves filtered colimits.
If $\C$ is (strongly) admissibly generated, then so is~$\D$.

\item \label{item--localization.admissibly.generated}
In the situation of \refth{symmetric.monoidal.localization}, suppose $\D$ is quasi-tractable.
Then the property of being (strongly) admissibly generated localizes.
\end{enumerate}
\xprop

\pf
\refit{admissibly.generated.weakly.sat}:
The proof is similar to the one of \refth{symmetric.weakly.saturated}.
In addition, we that an object~$X$ is small relative to some class~$\cell(T)$
if and only if it is small relative to its weak saturation \cite[Proposition 10.5.13]{Hirschhorn:Model}.

\refit{transfer.admissibly.generated}:
The cofibrant generation transfers to $\D$ by \refpr{transfer.basic}\refit{transfer.cofibrations}.
By Part~\refit{admissibly.generated.weakly.sat}, we only have to show that $\text{(co)dom}(F(I))$ are small with respect
to $\cell(Y \t_{\Sigma_n} s^{\pp n})$, where $s=F(t)$ are finite families of \emph{generating} cofibrations, i.e.,
$t$ are cofibrations in~$\C$.
By adjunction, this is equivalent to $\text{(co)dom}(I)$ being small with respect to
$$G(\cell(Y \t_{\Sigma_n} F(t)^{\pp n})) \subset \cell(G(Y \t_{\Sigma_n} F(t)^{\pp n})) = \cell (G(Y) \t_{\Sigma_n} t^{\pp n}),$$
which holds by assumption.

\refit{localization.admissibly.generated}: This is clear since $\cof_\C = \cof_\D$.
\xpf

By \refpr{strongly.admissibly.generated}, flatness and (symmetric) h-monoidality of $\Top$ only needs to be checked for generating cofibrations, which is easy.
Hence, $\Top$ is flat and symmetric h-monoidal (but not symmetric flat).

\subsection{Symmetric spectra}

The positive stable model structure on symmetric spectra valued in an abstract model category $\C$ is both symmetric flat and symmetric h-monoidal.
With a careful choice of the model structure on symmetric sequences, it is also symmetroidal.
As a special case, this shows that any model category is Quillen equivalent to one which is symmetric flat and symmetroidal.
For this, only mild conditions on $\C$ are necessary (such as flatness and h-monoidality, but not their symmetric counterparts).
See~\csp{\refsect{symmetricity.2}} for the precise statement.

\subsection{Equivariant homotopy theory}

For a symmetric monoidal model category $\C$, and a finite discrete group $G$,
we can work with $G$-equivariant homotopy theory in~$\C$ in two ways, namely by using the category $G\C$ of $G$-objects~(\refsect{finite.group.action}) in~$\C$
and the orbit presheaf category $\Orb_G \C$ defined as $\Fun(\Orb_G^\op, \C)$.
In this section, we list the (symmetric) monoidal properties of these two categories.

Recall from \cite{Stephan:Equivariant} that the objects of the {\it orbit category\/} $\Orb_G$ are $G$-orbits of the form $G/H$.
Here and below $H$ and~$K$ denote subgroups of $G$.
Morphisms are $G$-equivariant maps.
The obvious functors $\coprod_{H\subset G}\{*\}\to\Orb_G\gets\{G\}$ induce the evaluation functors in the following pair of adjunctions:
$$\xymatrix{\prod_H \C \ar@<.5ex>[rr] & & \Orb_G \C \ar@<.5ex>[ll]^{\ev} \ar@<.5ex>[rr]^{\ev_G} & & G\C. \ar@<.5ex>[ll].}$$

The left adjoint of $\ev$ is the unique cocontinuous functor
that sends an object $X \in \C$ concentrated in degree~$H$ for some fixed $H \subset G$ to $X \cdot \Hom_{\Orb_G}(-, G/H)$.
This object is denoted by $h_{G/H, X}$ and (by a slight abuse of language) called a \emph{representable object}.
The right adjoint of $\ev_G$ sends $X$ to the orbit object $G/K \mapsto X^K$,
where the superscript denotes the fixed point functor (i.e., the composition $G \C \r H \C \r \C$ of the restriction and limit functors).
There is a unique closed monoidal structure on $\Orb_G\C$ that satisfies
$$h_{G/H, X} \t h_{G/K, Y} = \coprod h_{G/H_i, X \t Y}, \eqlabel{orbit.products}$$
where $G/H \x G/K = \bigsqcup G/H_i$ is a decomposition of the cartesian (regarded as a $G$-space with the diagonal $G$-action) product into $G$-orbits.

%
The monoidal structure on $G\C$ is given by the one of $\C$ on the underlying level and the diagonal $G$-action.

\defi
The {\it projective model structure\/} on $\Orb_G \C$ is the one transferred from $\prod \C$ along this adjunction.
The {\it equivariant model structure\/} on $G\C$ is the one transferred from $\Orb_G \C$ (equivalently from $\prod \C$ along the composite adjunction).
That is, weak equivalences and (acyclic) fibrations on $G\C$
are maps $f$ such that $f^H$
is a weak equivalence, respectively, (acyclic) fibration for any $H \subset G$.
\xdefi


We now show that the orbit presheaf category always has very good monoidal properties.

\prop
\begin{enumerate}[(i)]
\item
The projective model structure on $\Orb_G \C$ exists whenever $\C$ is cofibrantly generated.
\item
\label{item--Orb.monoidal}
If $\C$ is combinatorial, quasi-tractable, pretty small, monoidal, or (symmetric) h-monoidal, then the same is true for $\Orb_G \C$.
If weak equivalences in $\C$ are stable under finite coproducts, the (symmetric) flatness of $\C$ implies the one of $\Orb_G \C$.
\end{enumerate}
\xprop

\pf
The first three properties of~\refit{Orb.monoidal} follow from \refpr{transfer.basic}.
The monoidality of $\Orb_G \C$ immediately follows from~\refeq{orbit.products} since the maps $h_{G/H, s}$, for a generating (acyclic) cofibration $s$ of $\C$, are the generating (acyclic) cofibrations of $\Orb_G \C$.

For h-monoidality, observe that for any $C \in \Orb_G \C$, we have $C\otimes(G/H\otimes s)=(C\otimes G/H)\otimes s$, which allows us to invoke the h-monoidality of~$\C$ objectwise.
For flatness, $y\pp(G/H\otimes i)=(y\otimes G/H)\pp i$, and $y \otimes G/H$ is an objectwise weak equivalence
because its components are coproducts of~$y$, thus they are weak equivalences by assumption.
The symmetric counterparts are treated similarly.
\xpf

By contrast, the category of $G$-objects requires much more restrictive hypotheses on $\C$ to behave well.
Below, the three technical cellularity conditions
are specifically designed to ensure the existence of a transferred model structure.
The conditions in the main statement below and part~\refit{equivariant.monoidal} are satisfied for $\C=\sSet$, $\Ch_\Q$, but not for $\sAb$ or $\Ch_\Z$.

\prop
\label{lemm--G.objects.model}
Suppose a model category $\C$ is finitely combinatorial,
the fixed point functor $(-)^H\colon G\C \r \C$ preserves pushouts along $G/K\otimes f$, where $H, K \subset G$, and $f$ is a generating cofibration in~$\C$,
and $Y \t (G/K)^H\to(Y \t G/K)^H$
is an isomorphism in~$\C$ for any $Y\in\C$, $K, H \subset G$.
\begin{enumerate}[(i)]
\item \label{item--equivariant.existence} \cite[Proposition~2.6]{Stephan:Equivariant}
The equivariant model structure on $G\C$ exists.
The adjunction $\Orb_G \C \rightleftarrows G\C$ is a Quillen equivalence.
\item
\label{item--equivariant.basic}
If $\C$ is combinatorial, quasi-tractable, pretty small, or (symmetric) h-monoidal, then so is $G\C$.
\item
\label{item--equivariant.monoidal}
Suppose that $\C$ is pretty small and that for any~$H$ the functor $(-)^H$ preserves pushouts along maps of the form $Y \t G/K \t s$, where $Y \in G\C$ and $s$ is a cofibration in $\C$.
If $\C$ is h-monoidal, then $G\C$ is h-monoidal.
If, moreover, $\C$ is flat, then $G\C$ is flat as well.
\end{enumerate}
\xprop

\pf
\refpr{transfer.basic} implies~\refit{equivariant.basic}.
By \refre{h.monoidal.weakly.sat}, h-monoidality only has to be checked for generating (acyclic) cofibrations in $G\C$.
These are of the form $G/K\otimes s$,
where $s$ is a generating (acyclic) cofibration of~$\C$.
For an arbitrary object $Y\in G\C$, the map $Y\otimes G/K\otimes s$ is an h-cofibration.
By \refle{i.cofibrations.detect}, it is enough to compute $H$-fixed points, which are (by assumption) $(Y\otimes G/K)^H\otimes s$.
This is an (acyclic) h-cofibration by the h-monoidality of~$\C$.
Similarly, flatness only has to be checked on generating cofibrations, where it follows by using the assumptions in the same way.
\xpf


\numberwithin{equation}{section}

\section{Applications to D-modules}

In this section, we establish model structures on operadic algebras in $\DD$-modules.
In particular, we recover the main theorem of \cite{diBrinoPistaloPoncin:Model}.

Let $X$ be a smooth algebraic variety over a field of characteristic zero.
(Entirely analogously, $X$ could be a complex manifold.) Let $\DD$ be the sheaf of differential operators on $X$.
The category $\Mod_\DD(X)$ is defined to be the category of left $\DD$-modules in the category $\PSh (X) := \PSh(X, \Ch(\Mod_\Q))$, i.e., presheaves of unbounded complexes of $\Q$-vector spaces.
The usual tensor product \cite[Proposition 1.2.9]{HottaTakeuchiTanisaki:D-modules} of left $\DD$-modules is denoted by $\t_\cO$.
We don't impose any quasi-coherence condition on $\DD$-modules.

The \emph{local projective model structure} on $\Mod_\DD(X)$ is constructed by means of the following adjunctions:
$$\prod_{U \subset X} \Ch(\Mod_\Q) \rightleftarrows \PSh(X)^{\proj} \rightleftarrows \PSh(X)^{\loc, \proj} \rightleftarrows \Mod_\cO(X) \rightleftarrows \Mod_\DD(X).\eqlabel{adjunction.Mod.D}$$
The projective model structure on $\Ch(\Mod_\Q)$ satisfies the following conditions~(\refsect{chain.complexes}):
\begin{center}(*)\qquad pretty small, tractable, symmetric h-monoidal, and symmetric flat.\end{center}
The first adjunction is the chain complex analogue of~\refeq{transfer.presheaf}.

\lemm
The projective model structure $\PSh(X)^\proj$ satisfies the properties (*).
\xlemm

\pf
Pretty smallness and tractability follow from \refpr{transfer.basic}.
\refth{symmetric.weakly.saturated} implies that it is enough to show the symmetric h-monoidality for \emph{generating} (acyclic) cofibrations $s$.
These are of the form $s = \Q[U] \t f$, where $U \subset X$ is open,
$\Q[-]\colon \PSh(\Set) \to \PSh(\Mod_\Q)$ denotes the $\Q$-linearization functor,
and $f \in \Ch(\Mod_\Q)$ is a generating (acyclic) cofibration.
The functor $\Q[-]$ is strong monoidal.
Therefore, for any presheaf of complexes~$Y$, any finite family of generating (acyclic) cofibrations $f$ of $\Ch(\Mod_\Q)$, and any finite family $U$ of open subsets of $X$, we have
$Y \t_{\Sigma_n} (\Q[U] \t f)^{\pp n} =
Y \t_{\Sigma_n} \left (\Q[U^{\x n}] \t f^{\pp n} \right )$.
(Here we are using the same notation as in \refde{multi}, so $U^{\x n}$ stands for $\prod_{i=1}^e U_i^{\x n_i}$.)
The evaluation functor is strong monoidal: for any open $V \subset X$, we get
$Y(V) \t_{\Sigma_n} \left ((U^{\x n})(V) \t f^{\pp n} \right ) = (Y(V) \t U^{\x n}(V)) \t_{\Sigma_n} f^{\pp n}$, which is an (acyclic) h-cofibration by symmetric h-monoidality of $\Ch(\Mod_\Q)$.
The symmetric flatness is shown similarly, using that $U^{\x n}(V)$ is either $\Q$ or $0$.
\xpf

The second adjunction in~\refeq{adjunction.Mod.D} is the monoidal Bousfield localization with respect to the Zariski topology, similarly to~\refsect{simplicial.presheaves}.
By \refpr{monoidal.localization} and \refth{symmetric.monoidal.localization}, the properties (*) remain valid.
The right adjoint in the third adjunction is the forgetful functor, along which the model structure is transferred to $\Mod_\cO(X)$.
This adjunction satisfies the conditions of \refth{commutative.monoid}, so the properties (*) transfer to $\cO$-modules.
The final right adjoint is again the forgetful functor.
Its left adjoint, $M \mapsto \DD \t_\cO M$, is not strong monoidal,
so we cannot apply \refth{commutative.monoid} as is, which forces us to give the following proof.

\theo
\label{theo--Mod.D.properties}
The local projective model structure on $\Mod_\DD(X)$ satisfies the properties (*).
\xtheo

\pf
Again invoking \refpr{transfer.basic} and \refth{symmetric.weakly.saturated}, the two symmetricity properties only have to be checked for generating (acyclic) cofibrations.
These are of the form $s = \DD \t_\cO f$, where $f$ is a generating (acyclic) cofibration of $\Mod_\cO(X)$.
For symmetric h-monoidality, we have to check that $Y \t_{\Ax_n, \cO} s^{\pp_\cO n}$
is an (acyclic) h-cofibration
for any multi-index $n$, any familiy $f$ of generating (acyclic) cofibrations and any $\Ax_n$-equivariant object $Y \in \Mod_\DD(X)$.
Here $\t_\cO$ and $\pp_\cO$ refer to the tensor product and pushout product of left $\DD$-modules.
Being a weak equivalence, and therefore also being an (acyclic) h-cofibration, is a local condition on $X$.
We may therefore assume that $\Theta$ is a free $\cO$-module.
Hence, $\DD$ is (non-canonically) isomorphic to a free $\cO$-module as well.
The forgetful functor $\Mod_\DD(X) \to \Mod_\cO(X)$ preserves pushouts and creates weak equivalences.
It therefore also detects (acyclic) h-cofibrations.
Consequently, it is enough to show that the above map, i.e.,
$Y \t_{\Sigma_n,\cO} (\DD \t_\cO f)^{\pp_\cO n}$, is an (acyclic) h-cofibration of $\cO$-modules.
Indeed, $\DD \t_\cO f$ is an (acyclic) cofibration of $\cO$-modules by the freeness of $\DD$ as an $\cO$-module, so we are done by the symmetric h-monoidality of $\Mod_\cO(X)$.

The same technique works for symmetric flatness.
\xpf

\rema
In order to obtain a similar statement for symmetroidality of $\DD$-modules, one has to work with a model structure whose cofibrations are local, i.e., satisfy descent.
We leave this to the reader.
\xrema

The following theorem extends the existence of the model structure on commutative $\DD$-algebras due to di~Brino, Pistalo, and Poncin \cite{diBrinoPistaloPoncin:Model}, which amounts to the case $O = \Comm$, the commutative operad.
It also provides a means for rectifying multiplicative structures on $\DD$-modules.

\coro
For any symmetric operad $O$ in $\Mod_\DD(X)$, the category of $O$-\hskip0pt algebras admits a model structure transferred along the forgetful functor
$$\Mod_\DD(X) \gets \Alg_O(\Mod_\DD(X)).$$
Moreover, any weak equivalence $O \to P$ of symmetric operads induces a Quillen equivalence
$$\Alg_O(\Mod_\DD(X)) \rightleftarrows \Alg_P(\Mod_\DD(X)).$$
\xcoro

\pf
By \refth{Mod.D.properties}, we can apply \cop{\refth{O.Alg}, \refth{O.strongly.admissible}, \refth{rect.colored.operad}}.
\xpf

\raggedright\rightskip0em plus \maxdimen

\bibliographystyle{dp}
\def\ZM#1{\href{https://zbmath.org/?q=an:#1}{Zbl #1}}
\def\MR#1{\href{http://www.ams.org/mathscinet-getitem?mr=#1}{MR#1}}
\bibliography{bib}

\end{document}